\documentclass[12pt]{amsart}
\usepackage{graphicx}
\usepackage{amssymb}
\usepackage{amsfonts}
\usepackage{amsmath,amscd}
\usepackage{array}
\usepackage{rotating}
\usepackage{epsfig}
\usepackage[matrix,frame]{xypic}

\newtheorem{example}{Example}[section]

\newtheorem{theorem}[example]{Theorem}

\newtheorem{corollary}[example]{Corollary}

\newtheorem{proposition}[example]{Proposition}

\newtheorem{lemma}[example]{Lemma}

\def\Proof{\noindent \it Proof -- \rm}
\def\qed{\hspace{3.5mm} \hfill \vbox{\hrule height 3pt depth 2 pt width 2mm}
\bigskip}

\def\Sym{{\bf Sym}}            % NCSF
\def\NCSF{{\bf Sym}}           % NCSF
\def\FQSym{{\bf FQSym}}        % permutations
\def\PBT{{\bf PBT}}            % arbres binaires planaires

\def\WQSym{{\bf WQSym}}        % Mots initiaux 
\def\maj{{\rm imaj\,}}
\def\maj{{\rm maj\,}}

\def\P{{\bf P}}
\def\M{{\bf M}}
\def\pack{{\rm pack}}
\def\PW{{\rm PW}}
\def\ev{{\rm Ev}}       % evaluation
\def\Std{{\rm Std}}     % standardisation
\def\cstd{{\rm cstd}}   % standardisation cyclique
\def\<{\langle}
\def\>{\rangle}
\def\NN{{\mathbb N}}    % entiers naturels
\def\ZZ{{\mathbb Z}}    % entiers relatifs
\def\gl{{\mathfrak gl}}

\def\F{{\bf F}}         % F de FQSym
\def\S{{\bf S}}         % S de NCSF
\def\G{{\bf G}}         % G de FQSym^*
\def\SG{{\mathfrak S}}  % groupe symetrique

\def\K{{\mathbb K}}
\def\X{{\bf X}}

\def\Des{\operatorname{Des}}
\def\maj{{\rm maj\,}}

\def\D{{\rm D}}
\def\Id{\operatorname{Id}}

\def\shuff#1#2{\mathbin{
\hbox{\vbox{ \hbox{\vrule \hskip#2 \vrule height#1 width 0pt
}%
\hrule}%
\vbox{ \hbox{\vrule \hskip#2 \vrule height#1 width 0pt
\vrule }%
\hrule}%
}}}

\long\def\psboxit#1#2{%
\begingroup\setbox0=\hbox{#2}%
\dimen0=\ht0 \advance\dimen0 by \dp0%
    \hbox{%
    \copy0%
    }%hbox
\endgroup%
}%psboxit
\def\SetTableau#1#2#3#4{%
  \gdef\Tabvrule{\vrule\vrule width-0.4pt}
  \gdef\Tabhrule{\hrule\hrule height-0.4pt}  
  \gdef\Tabstrut{\vrule height#1 depth#2 width0pt\relax}
  \gdef\Tabbox##1{\hbox to #3{\hskip0.4pt\hfill\Tabstrut$#4##1$\hfill}}
} %setTableau
\def\PetitTableau{\SetTableau{1.65ex}{0.55ex}{2.2ex}{\scriptstyle}}

\def\Case#1{\vcenter{\Tabhrule%
                   \hbox{\Tabvrule\Tabbox{#1}\Tabvrule}\Tabhrule}}
\def\GenTab#1{\vcenter{\halign{&$\Case{##}$\cr#1}}\egroup}
\def\Tableau{%         
  \bgroup%
  \let\ =\omit%
  \let\\=\cr%
  \offinterlineskip\GenTab}

\def\shuf{{\mathchoice{\shuff{7pt}{3.5pt}}%
{\shuff{6pt}{3pt}}%
{\shuff{4pt}{2pt}}%
{\shuff{3pt}{1.5pt}}}}%
\def\shuffle{\,\shuf\,}

\def\sev{{\rm sev}} % Evaluation signee
\def\BW{{\mathcal BW}} % Algebre de Hopf des mots signes tasses reguliers

\def\TT{{\mathcal T}} % Arbre plan associe a un mot
\def\ol#1{{\overline #1}}
\def\Sp{{SP}}
\def\Sm{{SM}}
\def\moinsu{{m(\epsilon)}}
\def\moinsup{{m(\epsilon')}}
\def\gf#1#2{\genfrac{}{}{0pt}{}{#1}{#2}}
\def\IN{{\rm Int}}
\def\Xb{{\bar X}}
\def\Ab{{\bar A}}

\def\gaudend{\prec}      % gauche dendriforme
\def\droitdend{\succ}    % droite dendriforme
\def\gautrid{\!\prec\!}   % gauche trigebre dendriforme
\def\miltrid{\circ}       % milieu trigebre dendriforme
\def\droittrid{\!\succ\!} % droite trigebre dendriforme
\def\PW{\rm PW} % mots tasses
\def\Nmax{{\rm NbMax}} % nb de max dans un mots tasses
\def\std{{\rm std}}
\def\cstd{{\bf Std}}
\def\AA{{\bf A}}
\def\alph{{\rm alph}}
\def\C{C}
\def\X{{\mathbb X}}

\def\NW{{\bf N}} % base duale des M_u dans WQSym
\def\NN{{\mathbb N}}

\def\evt{{\bf Ipack}}
\def\finerW{\geq}

\def\TD{{\mathfrak{TD}}}% Trigebre dendriforme libre sur un generateur
\def\MM{{\mathcal M}}  % Base standard de TD.

\def\spack{{\rm spack}}
\def\finer{\geq}

\title[Superization and $(q,t)$-specialization in combinatorial Hopf algebras]%
{Superization and $(q,t)$-specialization\\ in combinatorial Hopf algebras}

\author[J.-C.~Novelli and J.-Y.~Thibon]%
{Jean-Christophe Novelli and Jean-Yves Thibon}

\address[Novelli and Thibon]{Universit\'e Paris-Est \\Institut Gaspard Monge\\
5, Boulevard Descartes \\Champs-sur-Marne \\77454 Marne-la-Vall\'ee cedex 2 \\
FRANCE}
\email[Jean-Christophe Novelli]{novelli@univ-mlv.fr}
\email[Jean-Yves Thibon]{jyt@univ-mlv.fr}

\date{\today}

\begin{document}

\begin{abstract}
We extend a classical construction on symmetric functions, the superization
process, to several combinatorial Hopf algebras, and obtain analogs of the
hook-content formula for the $(q,t)$-specializations of various bases.
Exploiting the dendriform structures
yields in particular $(q,t)$-analogs of the Bj\"orner-Wachs
$q$-hook-length formulas for binary trees, and similar formulas for
plane trees.
\end{abstract}

\maketitle

%%%%%%%%%%%%%%%%%%%%%%%%%%%%%%%%%%%%%%%%%%%%%%%%%%%%%%%%%%%%%%%%%%%%%%%%%%%%%%%
%%%%%%%%%%%%%%%%%%%%%%%%%%%%%%%%%%%%%%%%%%%%%%%%%%%%%%%%%%%%%%%%%%%%%%%%%%%%%%%
%%%%%%%%%%%%%%%%%%%%%%%%%%%%%%%%%%%%%%%%%%%%%%%%%%%%%%%%%%%%%%%%%%%%%%%%%%%%%%%
\section{Introduction}

Combinatorial Hopf algebras are special graded and connected Hopf algebras
based on certain classes of combinatorial objects. 
There is no general agreement of what their precise definition should be,
but looking at their structure as well as to their existing applications,
it is pretty clear that they are to be regarded as generalizations of the
Hopf algebra of symmetric functions.

It is well known that one can define symmetric functions $f(X-Y)$
of a formal difference of alphabets. This can be interpreted either as the image
of the difference $\sum_i x_i -\sum_j y_j$ by the operator $f$
in the  $\lambda$-ring
generated by $X$ and $Y$, or, in Hopf-algebraic terms, as
$(\Id\otimes\tilde\omega)\circ\Delta (f)$, where $\Delta$ is the
coproduct and $\tilde\omega$ the antipode. And in slightly less pedantic
terms, this just amounts to replacing the power-sums $p_n(X)$ by
$p_n(X)-p_n(Y)$, a process already discussed at length in
Littlewood's book \cite[p. 100]{Lit}.

As is well known, the Schur functions $s_\lambda(X)$ are the characters
of the irreducible tensor representations of the general Lie algebra
$\gl(n)$. Similarly, the $s_\lambda(X-Y)$ are the characters of
the irreducible tensor representations of the general Lie superalgebras
$\gl(m|n)$ \cite{BeRe}. These symmetric functions are not positive sums of
monomials, and for this reason, one often prefers to use as characters
the so-called {\em supersymmetric functions} $s_\lambda(X|Y)$,
which are defined by $p_n(X|Y)=p_n(X)+(-1)^{n-1}p_n(Y)$ (see
\cite{Stem}), and are indeed positive sums of monomials: their
complete homogeneous functions are given by
\begin{equation}
\sigma_t(X|Y)=\sum_{n\ge 0}h_n(X|Y)t^n=\lambda_t(Y)\sigma_t(X)
=\prod_{i,j}\frac {1+ty_j}{1-tx_i}\,.
\end{equation}

Another (not unrelated) classical result on Schur functions is the
{\em hook-content} formula \cite[I.3 Ex. 3]{Mcd}, which gives in closed form
the specialization of a Schur function at the virtual alphabet
\begin{equation}
\frac{1-t}{1-q}=\frac1{1-q}-t\frac1{1-q}=1+q+q^2+\cdots -(t+tq+tq^2+\cdots)\,.
\end{equation}
This specialization was first considered by Littlewood \cite[Ch. VII]{Lit},
who obtained a factorized form for the result, but with possible
simplifications. The improved version known as the hook-content formula
\begin{equation}
s_\lambda\left(\frac{1-t}{1-q}\right)
=q^{n(\lambda)} \prod_{x\in\lambda}\frac{1-tq^{c(x)}}{1-q^{h(x)}}\,,
\end{equation}
which is a $(q,t)$-analog of the famous hook-length formula of
Frame-Robinson-Thrall \cite{FRT}, is due to Stanley \cite{Stan}.

The first example of a combinatorial Hopf algebra generalizing
symmetric functions is Gessel's algebra of quasi-symmetric functions
\cite{Ges}. Its Hopf structure was further worked out in \cite{MR,NCSF1},
and later used in \cite{NCSF2}, where two \emph{different} analogs
of the hook-content formula for quasi-symmetric functions are given.
Indeed, the notation
\begin{equation}
F_I\left(\frac{1-t}{1-q}\right)
\end{equation}
is ambiguous. It can mean (at least) two different things:
$$
\text{either}\ F_I\big(\, \frac{1}{1-q} \, \hat{\times}\, (1-t)\, \big)
\quad \hbox{or} \quad
F_I\big(\, (1-t) \,\hat{\times}\, \frac{1}{1-q}\, \big) \,,
$$
where $\hat{\times}$ denotes the ordered product of alphabets.
The second one is of the form $F_I(X-Y)$ (in the sense of \cite{NCSF2}),
but the first one is not (cf. \cite{NCSF2}).

In this article, we shall extend the notion of superization to several
combinatorial Hopf algebras. We shall start with $\FQSym$ (Free
quasi-symmetric functions, based on permutations), and our first result
(Theorem~\ref{GAAb}) will allow us to give new expressions and combinatorial
proofs of the $(q,t)$-specializations of quasi-symmetric functions.
Next, we extend these results to $\PBT$, the Loday-Ronco algebra of planar
binary trees, and obtain a $(q,t)$-analog of the Knuth and Bj\"orner-Wachs
hook-length formulas for binary trees. 
These results rely on the dendriform structure of $\PBT$. 
Exploiting in a similar way the tridendriform structure of $\WQSym$ (Word
quasi-symmetric functions, based on packed words, or set compositions), we
arrive at a $(q,t)$ analog of the formula of \cite{HNT2} counting packed words
according to the shape of their plane tree.

%%%%%%%%%%%%%%%%%%%%%%%%%%%%%%%%%%%%%%%%%%%%%%%%%%%%%%%%%%%%%%%%%%%%%%%%%%%%%%%
{\footnotesize {\it Acknowledgments.-}
This work has been partially supported by Agence Nationale de la Recherche,
grant ANR-06-BLAN-0380.
The authors would also like to thank the contributors of the MuPAD project,
and especially those of the combinat package, for providing the development
environment for this research (see~\cite{HT} for an introduction to
MuPAD-Combinat).
}

%%%%%%%%%%%%%%%%%%%%%%%%%%%%%%%%%%%%%%%%%%%%%%%%%%%%%%%%%%%%%%%%%%%%%%%%%%%%%%%
%%%%%%%%%%%%%%%%%%%%%%%%%%%%%%%%%%%%%%%%%%%%%%%%%%%%%%%%%%%%%%%%%%%%%%%%%%%%%%%
%%%%%%%%%%%%%%%%%%%%%%%%%%%%%%%%%%%%%%%%%%%%%%%%%%%%%%%%%%%%%%%%%%%%%%%%%%%%%%%
\section{Background}

%%%%%%%%%%%%%%%%%%%%%%%%%%%%%%%%%%%%%%%%%%%%%%%%%%%%%%%%%%%%%%%%%%%%%%%%%%%%%%%
\subsection{Noncommutative symmetric functions}

The reader is referred to \cite{NCSF1} for the basic theory
of noncommutative symmetric functions. The encoding of Hopf-algebraic
operations by means of sums, differences, and products of virtual
alphabets is fully explained in \cite{NCSF2}.

It is customary to reserve the letters $A$, $B,\dots$ for
noncommutative alphabets, and $X$, $Y,\dots$ for commutative
ones.

We recall that the Hopf algebra of noncommutative symmetric
functions is denoted by $\Sym$, or by $\Sym(A)$ if we consider the realization
in terms of an auxiliary alphabet. Bases of $\Sym_n$ are labelled by
compositions $I$ of $n$. The noncommutative complete and elementary functions
are denoted by $S_n$ and $\Lambda_n$, and the notation $S^I$ means
$S_{i_1}\cdots S_{i_r}$. The ribbon basis is denoted by $R_I$.
The notation $I\vDash n$ means that $I$ is a composition of $n$.
The conjugate composition is denoted by $I^\sim$.
The \emph{descent set} of $I$ is
$\Des(I) = \{ i_1,\ i_1+i_2, \ldots , i_1+\dots+i_{r-1}\}$.
The \emph{descent composition} of a permutation $\sigma\in\SG_n$
is the composition $I=D(\sigma)$ of $n$ whose descent set is the descent set
of $\sigma$.

The graded dual of $\Sym$ is $QSym$ (quasi-symmetric functions).
The dual basis of $(S^I)$ is $(M_I)$ (monomial), and that of $(R_I)$
is $(F_I)$.

The {\em evaluation} $\ev(w)$ of a word $w$ over a totally ordered alphabet
$A$ is the sequence $(|w|_a)_{a\in A}$ where $|w|_a$ is the number of
occurences of $a$ in $w$. The {\em packed evaluation} $I=\evt(w)$ is the
composition obtained by removing the zeros in $\ev(w)$.

The Hopf structures on $\Sym$ and $QSym$ allows one to mimic, up
to a certain extent, the $\lambda$-ring notation.
If $A$ is a totally ordered alphabet,
\begin{equation}
\sigma_t((1-q){A}) := \lambda_{-qt}(A) \sigma_t(A),
\end{equation}
\begin{equation}
\sigma_t \left(\frac{A}{1-q}\right)
 := \cdots \sigma_{q^2t}(A)\sigma_{qt}(A)\sigma_t(A)\,.
\end{equation}
We usually consider that our auxiliary variable
$t$ is of rank one, so that $\sigma_t(A)=\sigma_1(tA)$.

%%%%%%%%%%%%%%%%%%%%%%%%%%%%%%%%%%%%%%%%%%%%%%%%%%%%%%%%%%%%%%%%%%%%%%%%%%%%%%%
\subsection{Free quasi-symmetric functions}

The \emph{standardized word} $\std(w)$ of a word $w\in A^*$ is the permutation
obtained by iteratively scanning $w$ from left to right, and labelling
$1,2,\ldots$ the occurrences of its smallest letter, then numbering the
occurrences of the next one, and so on.
For example, $\std(bbacab)=341625$.
For a word $w$ on the alphabet $\{1,2,\ldots\}$, we denote by $w[k]$ the word
obtained by replacing each letter $i$ by the integer $i+k$.

Recall from \cite{NCSF6} that for an infinite totally ordered alphabet $A$,
$\FQSym(A)$ is the subalgebra of $\K\<A\>$ spanned by the polynomials
\begin{equation}
\G_\sigma(A)=\sum_{\std(w)=\sigma}w
\end{equation}
the sum of all words in $A^n$ whose standardization is the permutation
$\sigma\in\SG_n$. The multiplication rule is, for $\alpha\in\SG_k$ and
$\beta\in\SG_l$,
\begin{equation}\label{multG}
\G_\alpha \G_\beta
=
\sum_{\genfrac{}{}{0pt}{}{\gamma\in\SG_{k+l};\,\gamma=u\cdot v}
{\std(u)=\alpha,\std(v)=\beta}}\G_\gamma\,.
\end{equation}

The noncommutative ribbon Schur function $R_I\in\Sym$ is then
\begin{equation}
R_I=\sum_{\D(\sigma)=I}\G_\sigma\,.
\end{equation}
This defines a Hopf embedding $\Sym\rightarrow\FQSym$.
As a Hopf algebra, $\FQSym$ is self-dual. The scalar
product materializing this duality is the one for which
$(\G_\sigma\,,\,\G_\tau)=\delta_{\sigma,\tau^{-1}}$
(Kronecker symbol).
Hence, $\F_\sigma:=\G_{\sigma^{-1}}$ is the dual basis
of $\G$. 

The {\em internal product} $*$ of $\FQSym$ is induced by
composition $\circ$ in $\SG_n$ in the basis~$\F$, that is,
\begin{equation}
\F_\sigma * \F_\tau = \F_{\sigma\circ\tau}\quad\text{and}\quad
\G_\sigma * \G_\tau = \G_{\tau\circ\sigma}\,.
\end{equation}
Each subspace $\Sym_n$ is stable under this operation, and anti-isomorphic to
the descent algebra $\Sigma_n$ of $\SG_n$.

The transpose of the Hopf embedding $\Sym\rightarrow\FQSym$
is the {\em commutative image} $\F_\sigma \mapsto \F_\sigma(X)=F_I(X)$,
where $I$ is the descent composition of $\sigma$,	
and $F_I$ is Gessel's fundamental basis of $QSym$.

%%%%%%%%%%%%%%%%%%%%%%%%%%%%%%%%%%%%%%%%%%%%%%%%%%%%%%%%%%%%%%%%%%%%%%%%%%%%%%%
%%%%%%%%%%%%%%%%%%%%%%%%%%%%%%%%%%%%%%%%%%%%%%%%%%%%%%%%%%%%%%%%%%%%%%%%%%%%%%%
%%%%%%%%%%%%%%%%%%%%%%%%%%%%%%%%%%%%%%%%%%%%%%%%%%%%%%%%%%%%%%%%%%%%%%%%%%%%%%%
\section{Free super-quasi-symmetric functions}

%%%%%%%%%%%%%%%%%%%%%%%%%%%%%%%%%%%%%%%%%%%%%%%%%%%%%%%%%%%%%%%%%%%%%%%%%%%%%%%
\subsection{Supersymmetric functions}

As already mentioned in the introduction, in the $\lambda$-ring notation, the
definition of \emph{supersymmetric functions} is transparent. If $X$ and $\Xb$
are two independent infinite alphabets, the \emph{superization} $f^\#$ of
$f\in\Sym$ is
\begin{equation}
f^\# := f(X\,|\,\Xb) = f(X-q\Xb)|_{q=-1},
\end{equation}
where $f(X-q\Xb)$ is interpreted in the $\lambda$-ring sense
($p_n(X-q\Xb):=p_n(X)-q^np_n(\Xb)$), $q$ being of rank one, so that
$p_n(X|\Xb)=p_n(X)-(-1)^np_n(\Xb)$.
This can also be written as an internal product
\begin{equation}
f^\# := f*\sigma_1^\#,
\end{equation}
where $\sigma_1^\#=\sigma_1(X-q\Xb)|_{q=-1} = \lambda_1(\Xb)\sigma_1(X)$,
and the internal product is extended to the algebra generated by $Sym(X)$ and
$Sym(\Xb)$ by means of the splitting formula
\begin{equation}
\label{split}
(f_1\dots f_r)*g = \mu_r\cdot (f_1\otimes\dots\otimes f_r)*_r \Delta^r g,
\end{equation}
and the rules
\begin{equation}
\label{rules}
\sigma_1 * f = f * \sigma_1,\qquad
{\overline\sigma_1} * {\overline\sigma_1} = \sigma_1.
\end{equation}

%%%%%%%%%%%%%%%%%%%%%%%%%%%%%%%%%%%%%%%%%%%%%%%%%%%%%%%%%%%%%%%%%%%%%%%%%%%%%%%
\subsection{Noncommutative supersymmetric functions}

The same can be done with noncommutative symmetric functions.
We need two independent infinite totally ordered alphabets $A$ and $\Ab$ and
we define $\NCSF(A|\Ab)$ as the subalgebra of the free product
${\NCSF^{(2)}}:=\NCSF(A)\star\NCSF(\Ab)$ generated by $S_n^\#$ where
\begin{equation}
\sigma_1^\# = 
\bar\lambda_1\sigma_1
=
\sum_{I=(i_1,\ldots,i_r)} (-1)^{i_1+\dots+i_r-r} \overline{S^{i_1\dots i_r}}
              S^{i_{r+1}}.
\end{equation}
For example,
\begin{equation}
S_1^\# = S^1 + S^{\ol1},
\qquad
S_2^\# = S^2 + S^{{\ol1}1} - S^{\ol2} + S^{\ol1\ol1},
\end{equation}
\begin{equation}
S_3^\# = S^3 + S^{\ol12} + S^{\ol1\ol11} - S^{\ol21} + S^{\ol1\ol1\ol1}
         - S^{\ol2\ol1}  - S^{\ol1\ol2} + S^{\ol3}.
\end{equation}

We shall denote the generators of $\NCSF^{(2)}$ by
$S_{(i,\epsilon)}$ where $\epsilon=\{\pm1\}$, so that
$S_{(i,1)}=S_i$ and $S_{(i,-1)}= S_{\ol i}$.

The corresponding basis of $\NCSF^{(2)}$ is then written
\begin{equation}
S^{(I,\epsilon)}=
S^{(i_1,\dots,i_r),(\epsilon_1,\dots,\epsilon_r)} :=
S_{(i_1,\epsilon_1)} S_{(i_2,\epsilon_2)} \dots S_{(i_r,\epsilon_r)},
\end{equation}
where $I=(i_1,\dots,i_r)$ is a composition and
$\epsilon=(\epsilon_1,\dots,\epsilon_r)\in\{\pm 1\}^r$ is a vector of signs.

Again, we extend the internal product by formulas~(\ref{split})
and~(\ref{rules}) where, now, $f_1$, $\dots,$ $f_r$, $g\in\NCSF^{(2)}$, and
$\sigma_1=\sigma_1(A)$, $\overline{\sigma_1}=\sigma_1(\Ab)$.
The resulting algebra is isomorphic to the Mantaci-Reutenauer algebra of type
$B$~\cite{ManR}.
We define the superization of $f\in\NCSF$ by
\begin{equation}
\label{sncsf}
f^\# := f * \sigma_1^\# = f(A-q\Ab)|_{q=-1} = f * (\ol\lambda_1\sigma_1).
\end{equation}

%%%%%%%%%%%%%%%%%%%%%%%%%%%%%%%%%%%%%%%%%%%%%%%%%%%%%%%%%%%%%%%%%%%%%%%%%%%%%%%
\subsection{Super-quasi-symmetric functions}

There are two natural and nonequivalent choices for defining
super-quasi-symmetric functions. The first one is to set
$F(X|\Xb)=F(X-q\Xb)|_{q=-1}$ as in~\cite{HHL}.
The second one is obtained by commutative image from the free
super-quasi-symmetric functions to be defined below.
Let us note that super-quasi-symmetric functions have been recently
interpreted as characters of certain abstract crystals of the Lie
superalgebras $\gl(m|n)$ \cite{Kwon}.

%%%%%%%%%%%%%%%%%%%%%%%%%%%%%%%%%%%%%%%%%%%%%%%%%%%%%%%%%%%%%%%%%%%%%%%%%%%%%%%
\subsection{Free super-quasi-symmetric functions}

The expressions (\ref{sncsf}) are still well-defined for an arbitrary
$f\in\FQSym$.
We can define
$\FQSym(A|\Ab)$ as the subalgebra of the free product $\FQSym(A)\star\FQSym(\Ab)$
spanned  by
\begin{equation}
\G_\sigma^\# := \G_\sigma(A|\Ab)=\G_\sigma * \sigma_1^\#.
\end{equation}
Again, $*$ is extended to the free product by conditions (\ref{split})
(valid only if $g\in\Sym^{(2)}$, which is enough), and (\ref{rules}).
This free product can be interpreted as $\FQSym^{(2)}$, the algebra of free
quasi-symmetric functions of level 2, as defined in \cite{FQSl}.
Let us set
\begin{gather}
A^{(0)}=A = \{a_1<a_2<\ldots<a_n<\ldots\}\,,\\
A^{(1)}=\Ab = \{\ldots <\bar a_n <\ldots <\bar a_2<\bar a_1\}\,,
\end{gather}
order $\AA=\bar A\cup A$ by $\bar a_i<a_j$ for all $i,j$,
and denote by $\std$  the
standardization of signed words with respect to this order.
We also need the signed standardization $\cstd$, defined as follows.
Represent a signed word ${\bf w}\in\AA^n$ by a pair
$(w,\epsilon)$, where $w\in A^n$ is the underlying unsigned word, and
$\epsilon\in\{\pm 1\}^n$ is the vector of signs. Then
$\cstd(w,\epsilon)=(\std(w),\epsilon)$.

We denote by $\moinsu$  the number of entries $-1$ in $\epsilon$.

A basis of $\FQSym^{(2)}$ is given by
\begin{equation}
\G_{\sigma,\epsilon} := \sum_{\cstd({\bf w})=(\sigma,\epsilon)} {\bf w}
\quad\in\ZZ\langle\AA\rangle.
\end{equation}
and the internal product obtained from (\ref{split}-\ref{rules}) 
coincides with the one of \cite{FQSl},
so that it is in fact always well-defined. 
In particular, viewing signed permutations as elements of
the group $\{\pm 1\}\wr\SG_n$,
\begin{equation}
\G_{\alpha,\epsilon}* \G_{\beta,\eta}=\G_{(\beta,\eta)\circ(\alpha,\epsilon)}
=\G_{\beta\circ\alpha, (\eta\alpha)\cdot \epsilon}
\end{equation}
with $\eta\alpha=(\eta_{\alpha(1)},\ldots,\eta_{\alpha(n)})$
and $\epsilon\cdot \eta=(\epsilon_1\eta_1,\ldots,\epsilon_n\eta_n)$.
\begin{theorem}
\label{GAAb}
The expansion of $\G_\sigma(A|\Ab)$ on the basis $\G_{\tau,\epsilon}$ is
\begin{equation}\label{Gdiese}
\G_\sigma(A|\Ab)=
\sum_{\std(\tau,\epsilon)=\sigma}\G_{\tau,\epsilon}\,.
\end{equation}
\end{theorem}

\Proof This is clear for $\sigma=12\ldots n$:
\begin{equation}
\sum_n \G_{12\ldots n}(A|\Ab) = \bar\lambda_1\cdot\sigma_1
=\sum_{\gf{i_1<i_2<\ldots <i_k}{j_1\le j_2\le\ldots\le j_l}}
\bar a_{i_1}\bar a_{i_2}\cdots \bar a_{i_k}a_{j_1}a_{j_2}\cdots a_{j_l}\,,
\end{equation}
and writing 
\begin{equation}
\G_\sigma(A|\Ab)=\G_\sigma * (\bar\lambda_1\cdot\sigma_1)
=\sum_{\std(\tau,\epsilon)=12\cdots n}\G_{\tau\sigma,\epsilon\sigma}
=\sum_{\std(\tau,\epsilon)=\sigma}\G_{\tau,\epsilon}\,,
\end{equation}
we obtain (\ref{Gdiese}).
\qed

%%%%%%%%%%%%%%%%%%%%%%%%%%%%%%%%%%%%%%%%%%%%%%%%%%%%%%%%%%%%%%%%%%%%%%%%%%%%%%%
\subsection{The canonical projection}

We have an obvious projection
\begin{equation}
\FQSym(A|\Ab) \to \FQSym(A)
\end{equation}
consisting in setting $\Ab=A$. One can even describe the refined map
\begin{equation}
\eta_t(\G_\sigma^\#) = \G_\sigma(A|tA).
\end{equation}

\begin{corollary}
\label{GAtA}
In the special case $\Ab=tA$, one gets
\begin{equation}
\G_\sigma(A|tA) = \sum_{\std(\tau,\epsilon)=\sigma} t^{\moinsu} \G_\tau(A).
\end{equation}
\end{corollary}

\Proof
This follows from (\ref{Gdiese}).
\qed

\begin{example}
{\rm
We have
\begin{equation}
\G_{12}(A|tA) = (1+t) (\G_{12} + t\G_{21}),
\quad\quad
\G_{21}(A|tA) = (1+t) (\G_{21} + t\G_{12})
\end{equation}
\begin{equation}
\begin{split}
\G_{4132}(A|tA) = (1+t)
(&\G_{4132} + t\G_{3421} + t\G_{4231} + t \G_{4321} \\
&+ t^2\G_{2413}+t^2 \G_{3412} + t^2\G_{4312} + t^3\G_{1423}).\label{ex26}
\end{split}
\end{equation}
Indeed, (\ref{ex26}) is obtained from the 16 signed permutations whose
standardized word is $4132$: 
\begin{equation}
\begin{split}
&
4132,\ 4\ol132,\ \
3\ol421,\ 3\ol42\ol1,\ \
4\ol231,\ 4\ol23\ol1,\ \
4\ol321,\ 4\ol32\ol1,\\
&2413,\ 24\ol13,\ \
3\ol41\ol2,\ 3\ol4\ol1\ol2,\ \
4\ol31\ol2,\ 4\ol3\ol1\ol2,\ \
1\ol4\ol2\ol3,\ \ol1\ol4\ol2\ol3.
\end{split}
\end{equation}
}
\end{example}

Summing over a descent class, we obtain
\begin{corollary}
\label{cor2FAt}
\begin{equation}
 R_I(A|\Ab) = \sum_{\C(J,\epsilon)=I} R_{J,\epsilon},
\end{equation}
where $R_{J,\epsilon}$ is the signed ribbon Schur function defined as
in~\cite{FQSl} and $\C(J,\epsilon)$ is the composition whose descents are the
descents of any signed permutation $(\sigma,\epsilon)$ where $\sigma$ is of
shape $J$.
\end{corollary}

Substituting $\Ab=tA$ yields
\begin{equation}
R_I(A|tA) = \sum_{\C(J,\epsilon)=I} t^{\moinsu} R_J(A),
\end{equation}
which allows us to recover a formula of \cite{NCSF2} 
(in \cite{NCSF2}, the exponent $b(I,J)$ is incorrectly stated).
Recall that a \emph{peak} of a composition is a cell of its diagram having no
cell to its right or on its top and that a \emph{valley} is a cell having no
cell to its left or at its bottom.

\begin{corollary}[\cite{NCSF2}, (121)]
\begin{equation}
 R_I(A|tA) = \sum_{J} (1+t)^{v(J)} t^{b(I,J)} R_I(A),
\end{equation}
where the sum is over all compositions $I$ which have either a peak or a
valley at each peak of $J$.
The power of $1+t$ is given by the number of valleys $v(J)$ of $J$ and the
power  of $t$ is the number of descents of $J$ that are not descents
of $I$ plus the number of descents $d$ of $I$ such that neither $d$ nor $d-1$
are descents of $J$.
\end{corollary}

\Proof
This is best understood at the level of permutations. 
First, the coefficient of
$R_J(A)$ is equal, by definition,  to the number of signed
permutations of shape $I$ whose underlying (unsigned) permutation is of
shape $J$.
Now, on the ribbon diagram of a permutation of shape $J$, in order to obtain a signed
permutation of shape $I$, we distinguish three kinds of cells: those
which must have a plus sign, those which must have a minus sign, and
those which can have both signs.
The valleys of $J$ can get any sign without changing their final shape whereas
all other cells have a  fixed plus or  minus sign, depending on $I$ and
$J$, thus explaining the coefficient $(1+t)^{v(J)}$.
The cells which must have a minus sign are either the descents of
$J$ that are not descents of $I$ or the descents $d$ (plus one) of $I$ such
that neither $d$ nor $d-1$ are descents of $J$, whence the power of
$t$. Indeed, it is enough to determine the correct
power of $t$ in the middle cell for all pairs of compositions of  $3$, 
since it depends only on the relative positions of their adjacent cells in $I$ and
$J$.
\qed

%%%%%%%%%%%%%%%%%%%%%%%%%%%%%%%%%%%%%%%%%%%%%%%%%%%%%%%%%%%%%%%%%%%%%%%%%%%%%%%
\subsection{The dual transformation}

Corollary~\ref{GAtA} is equivalent, up to substituting $-t$ to $t$, to a
combinatorial description of
\begin{equation}
\G_\sigma((1-t)A)=\G_\sigma(A) * \sigma_1((1-t)A).
\end{equation}
Let $\eta_t^*$ be the adjoint of $\eta_t$.  We can consistently set
\begin{equation}
\F_\sigma(A\cdot(1-t)) := \eta_t^* (\F_\sigma(A)),
\end{equation}
since the noncommutative Cauchy formula reads
\begin{equation}
\begin{split}
\sigma_1(A\cdot (1-t)\cdot B)
&= \sum_{\alpha} \F_\alpha(A\cdot (1-t))\G_\alpha(B) \\
&= \sum_{\beta} \F_\beta(A)\G_\beta((1-t)B).
\end{split}
\end{equation}
Writing
\begin{equation}
\begin{split}
\G_\beta((1-t)B)
&= \G_\beta(B) *  S_n((1-t)B)\\
&=  \sum_{k=0}^{n-1} (1-t) (-t)^k \G_\beta * R_{(1^k,n-k)}\\
&= \sum_{k=0}^{n-1} (1-t) (-t)^k
   \sum_{\Des(\tau)=\{1,\dots,k\}} \G_{\tau\circ\beta}(B),
\end{split}
\end{equation}
we have, setting $\sigma=\tau\circ\beta$
\begin{equation}
\sigma_1(A\cdot (1-t)\cdot B)
= 1 + \sum_{|\sigma|\geq1}
\left(  \sum_{\Des(\tau)=\{1,\dots,k\}} (1-t) (-t)^k
        \F_{\tau^{-1}\circ\sigma}(A)
\right)
\G_\sigma(B),
\end{equation}
so that
\begin{equation}
\F_\sigma(A\cdot (1-t)) 
= \sum_{k=0}^{n-1} (1-t)(-t)^k
   \sum_{\Des(\tau)=\{1,\dots,k\}} \F_{\tau^{-1}\circ\beta}(A).
\end{equation}

\begin{theorem}
\label{FAt}
In terms of signed permutations, we have
\begin{equation}
 \F_\sigma(A\cdot (1-t))
 = \sum_{\epsilon\in\{\pm1\}^n} (-t)^{\moinsu} \F_{\std(\sigma,\epsilon)}(A).
\end{equation}
\end{theorem}

\Proof
This follows from the above discussion and Corollary~\ref{GAtA}.
\qed

From now on, we denote by $\X$ the alphabet
$\frac{|1-t}{1-q|}:=\frac1{1-q}\, \hat\times\, (1-t)$.

\begin{corollary}
\label{corFAt}
Specializing $A=\frac{1}{1-q}$, we obtain
\begin{equation}
\begin{split}
 \F_\sigma(\X)
 &= F_{\D(\sigma)}(\X)
    \text{\rm\ \ in the notation of~\cite{NCSF2}}\\
 &= \frac{1}{(q)_n}\sum_{\epsilon\in\{\pm1\}^n} (-t)^{\moinsu} q^{\maj(\sigma,\epsilon)},
\end{split}
\end{equation}
where $\maj(w)=\sum_{i\in \Des(w)} i$, and
$\Des(w)=\{i| w_i>w_{i+1}\}$.
\end{corollary}

%%%%%%%%%%%%%%%%%%%%%%%%%%%%%%%%%%%%%%%%%%%%%%%%%%%%%%%%%%%%%%%%%%%%%%%%%%%%%%%
\subsection{A hook-content formula in $\FQSym$}

Let us denote by $\Sp_i$ the set of words $\epsilon\in\{\pm1\}^n$ where
$\epsilon_i=1$
and by $\Sm_i$ the set of words $\epsilon\in\{\pm1\}^n$ where $\epsilon_i=-1$.

Let $\phi_i$ be the involution on signed permutations $(\sigma,\epsilon)$ which
changes the sign of $\epsilon_i$ and leaves the rest unchanged.

\begin{lemma}
\label{lemFsig}
Let $(\sigma,\epsilon)$ be a signed permutation such that $\epsilon_i=1$ and let
$(\sigma,\epsilon')=\phi_i(\sigma,\epsilon)$. Then
\begin{equation}
(-t)^{\moinsup}q^{maj(\sigma,\epsilon')} = (-t) \frac{q^{(i-1)x_i}}{q^{iy_i}}
(-t)^{\moinsu}q^{maj(\sigma,\epsilon)},
\end{equation}
where $x_i=0$ if $\sigma_{i-1}>\sigma_i$ and $x_i=1$ otherwise,
and $y_i=0$ if $\sigma_{i}<\sigma_{i+1}$ and $y_i=1$ otherwise.
By convention, $x_1=0$ and $y_n=0$, which is equivalent to fix
$\sigma_0=\sigma_{n+1}=+\infty$.
\end{lemma}

\Proof
The factor $(-t)$  is obvious.
The difference between the $q$-statistics of both words 
depends only on the descents at
position $i-1$ and position $i$.
Let us discuss  position $i-1$ (value of $x_i$).
If $\sigma_{i-1}>\sigma_i$, we have
\begin{equation}
-\sigma_{i-1}<-\sigma_{i}<\sigma_{i}<\sigma_{i-1},
\end{equation}
so that there is a descent at position $i-1$ in $(\sigma,\epsilon)$ iff there is a
descent at the same position in $(\sigma,\epsilon')$. This proves the case
$x_i=0$.

Now, if $\sigma_{i-1}>\sigma_i$, we have
\begin{equation}
-\sigma_{i}<-\sigma_{i-1}<\sigma_{i-1}<\sigma_{i},
\end{equation}
so that there is no descent at position $i-1$ in $(\sigma,\epsilon)$ and there is a
descent at the same position in $(\sigma,\epsilon')$. This proves the case $x_i=1$.
The discussion of position $i$ is similar.
\qed

\begin{theorem}
\label{HookFQSym}
Let $\sigma\in\SG_n$. Then
\begin{equation}
\label{EqFSig}
\F_\sigma(\X)
 = q^{\maj(\sigma)}
    \prod_{i=1}^n \frac{1-q^{(i-1)x_i-iy_i}\,t}{1-q^i}\\
 = \prod_{i=1}^n \frac{q^{i y_i} -q^{(i-1)x_i}\,t}{1-q^i},
\end{equation}
where $x_i$ and $y_i$ are as in Lemma~\ref{lemFsig}.
\end{theorem}
This gives an analog of the hook-content formula,
where the hook-length of cell $\# i$ is its ``ribbon length'' $i$,
and its ``content'' is $c_i=(i-1)x_i-iy_i$.

\Proof
Thanks to Lemma~\ref{lemFsig}, we have
\begin{equation}
(q)_n \F_\sigma(\X)
 = \sum_{\epsilon\in \{\pm1\}^n} (-t)^{\moinsu}q^{\maj(\sigma,\epsilon)}
 = \left(1-t\frac{q^{(i-1)x_i}}{q^{iy_i}} \right)
\sum_{\epsilon\in \Sp_i} (-t)^{\moinsu}q^{\maj(\sigma,\epsilon)},
\end{equation}
since each signed permutation $(\sigma,\epsilon')$ with $\epsilon'_i=-1$ gives
the same contribution as $\phi_i(\sigma,\epsilon')$ up to the factor
involving $x_i$ and $y_i$.
The same can be done for signed permutations such that $\epsilon'_i=1$ and
$\epsilon'_j=-1$, so that the whole expression factors and gives the first
formula of (\ref{EqFSig}). The second expression is clearly equivalent.
\qed

%%%%%%%%%%%%%%%%%%%%%%%%%%%%%%%%%%%%%%%%%%%%%%%%%%%%%%%%%%%%%%%%%%%%%%%%%%%%%%%
\subsection{Graphical representations}

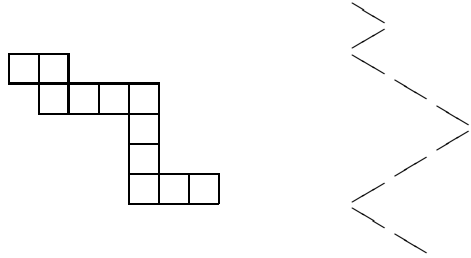
\begin{figure}[ht]
\entrymodifiers={+<4pt>}
\begin{equation*}
\PetitTableau
\PetitTableau
\Tableau{ & \\\ & & & &\\\ &\ &\ &\ & \\\ &\ &\ &\ & \\
\ &\ &\ &\ & & &\\}
\qquad
\qquad
\vcenter{\xymatrix@C=-7mm@R=+2mm{
*{} & *{} & *{} & {}\ar@{-}[dl] \\
*{} & *{} & {}\ar@{-}[dr] \\
*{} & *{} & *{} & {}\ar@{-}[dl]\\
*{} & *{} & {}\ar@{-}[dl] \\
*{} & {}\ar@{-}[dl] \\
{}\ar@{-}[dr] \\
*{} & {}\ar@{-}[dr] \\
*{} & *{} & {}\ar@{-}[dr] \\
*{} & *{} & *{} & {}\ar@{-}[dl] \\
*{} & *{} & {}\ar@{-}[dl] \\
*{} & {} \\
}}
\end{equation*}
\caption{\label{fig-compo-arbre}
The mirror shape of $\sigma=(5, 6, 7, 4, 3, 2, 8, 9, 10, 1, 11)$
and its representation as a binary tree.}
\end{figure}

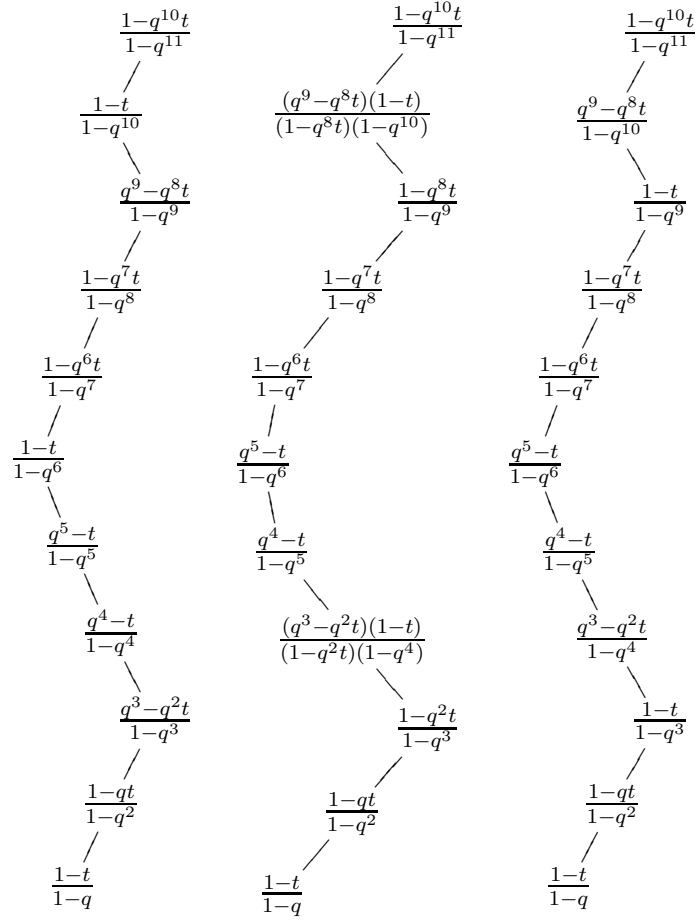
\begin{figure}[ht]
{
\entrymodifiers={+<4pt>}
\begin{equation*}
\vcenter{\xymatrix@C=-5mm@R=+4mm{
*{} & *{} & *{} & \frac{1-q^{10}t}{1-q^{11}}\ar@{-}[dl] \\
*{} & *{} & \frac{1-t}{1-q^{10}}\ar@{-}[dr] \\
*{} & *{} & *{} & \frac{q^9-q^8t}{1-q^9}\ar@{-}[dl]\\
*{} & *{} & \frac{1-q^7t}{1-q^8}\ar@{-}[dl] \\
*{} & \frac{1-q^6t}{1-q^7}\ar@{-}[dl] \\
\frac{1-t}{1-q^6}\ar@{-}[dr] \\
*{} & \frac{q^5-t}{1-q^5}\ar@{-}[dr] \\
*{} & *{} & \frac{q^4-t}{1-q^4}\ar@{-}[dr] \\
*{} & *{} & *{} & \frac{q^3-q^2t}{1-q^3}\ar@{-}[dl] \\
*{} & *{} & \frac{1-qt}{1-q^2}\ar@{-}[dl] \\
*{} & \frac{1-t}{1-q} \\
}}
\quad
\vcenter{\xymatrix@C=-7mm@R=+4mm{
*{} & *{} & *{} & \frac{1-q^{10}t}{1-q^{11}}\ar@{-}[dl] \\
*{} & *{} & \frac{(q^9-q^8t)(1-t)}{(1-q^8t)(1-q^{10})}\ar@{-}[dr] \\
*{} & *{} & *{} & \frac{1-q^8t}{1-q^9}\ar@{-}[dl]\\
*{} & *{} & \frac{1-q^7t}{1-q^8}\ar@{-}[dl] \\
*{} & \frac{1-q^6t}{1-q^7}\ar@{-}[dl] \\
\frac{q^5-t}{1-q^6}\ar@{-}[dr] \\
*{} & \frac{q^4-t}{1-q^5}\ar@{-}[dr] \\
*{} & *{} & \frac{(q^3-q^2t)(1-t)}{(1-q^2t)(1-q^4)}\ar@{-}[dr] \\
*{} & *{} & *{} & \frac{1-q^2t}{1-q^3}\ar@{-}[dl] \\
*{} & *{} & \frac{1-qt}{1-q^2}\ar@{-}[dl] \\
*{} & \frac{1-t}{1-q} \\
}}
\quad
\vcenter{\xymatrix@C=-5mm@R=+4mm{
*{} & *{} & *{} & \frac{1-q^{10}t}{1-q^{11}}\ar@{-}[dl] \\
*{} & *{} & \frac{q^9-q^8t}{1-q^{10}}\ar@{-}[dr] \\
*{} & *{} & *{} & \frac{1-t}{1-q^9}\ar@{-}[dl]\\
*{} & *{} & \frac{1-q^7t}{1-q^8}\ar@{-}[dl] \\
*{} & \frac{1-q^6t}{1-q^7}\ar@{-}[dl] \\
\frac{q^5-t}{1-q^6}\ar@{-}[dr] \\
*{} & \frac{q^4-t}{1-q^5}\ar@{-}[dr] \\
*{} & *{} & \frac{q^3-q^2t}{1-q^4}\ar@{-}[dr] \\
*{} & *{} & *{} & \frac{1-t}{1-q^3}\ar@{-}[dl] \\
*{} & *{} & \frac{1-qt}{1-q^2}\ar@{-}[dl] \\
*{} & \frac{1-t}{1-q} \\
}}
\end{equation*}
}
\caption{\label{figEqsFQS}The three hook-content formulas for the permutation
$(5, 6, 7, 4, 3, 2, 8, 9, 10, 1, 11)$: signed permutations (left diagram),
induction (middle diagram), and simplification of the induction
(right diagram).}
\end{figure}

We shall see later that (\ref{EqFSig}) is the special case of formula
(\ref{EqArb}) for binary trees, when the tree is a zig-zag line. This is why
we have chosen to represent graphically $\F_\sigma(\X)$ with hook-content type
factors in the following way:
let the \emph{mirror shape} of a permutation $\sigma$ be the mirror image of
its descent composition. We represent it as the binary tree in which each
internal node has only one subtree, depending on whether the corresponding cell
of the composition is followed by a cell to its right or to its bottom.
For example, with $\sigma=(5, 6, 7, 4, 3, 2, 8, 9, 10, 1, 11)$, 
the shape is $(3,1,1,4,4)$, the mirror shape is $(2,4,1,1,3)$ and its binary
tree is shown on Figure~\ref{fig-compo-arbre}.
Theorem \ref{HookFQSym} can be visualized by placing into the $i$th node (from
bottom to top) the $i$-th factor of $\F_\sigma(\X)$ in Equation~(\ref{EqFSig}).
For example, the first tree of Figure~\ref{figEqsFQS} shows the expansion of
$\F_\sigma(\X)$ with the hook-content factors of
$\sigma=(5, 6, 7, 4, 3, 2, 8, 9, 10, 1, 11)$.
We shall see  two alternative hook-content formulas for
$\F_\sigma(\X)$. The first one is obtained from an induction formula expressing
$\F_\sigma(\X)$ from
$\F_{\Std(\sigma_1\dots\sigma_{n-1})}(\X)$,
and follows directly from Theorem~\ref{HookFQSym}.

\begin{corollary}
\label{propFsig}
Let $\partial\F_\sigma(\X):=\F_{\Std(\sigma_1\dots\sigma_{n-1})}(\X)$
as in \cite{HNT2}. Then,
\begin{equation}
\label{FSigrecc}
\F_\sigma(\X) = 
\partial\F_\sigma(\X)\times
\left\{
\begin{array}{ll}
\displaystyle \frac{1-q^{n-1}t}{1-q^n}
& \text{\rm if $\sigma_{n-1}\!<\!\sigma_n$,} \\[.3cm]
\displaystyle\frac{q^{n-1}-t}{1-q^n}
& \text{\rm if $\sigma_{n-2}\!>\!\sigma_{n-1}\!>\!\sigma_n$}, \\[.3cm]
\displaystyle \frac{(q^{n-1}-q^{n-2}t)(1-t)}{(1-q^{n-2}t)(1-q^n)}
& \text{\rm if $\sigma_{n-2}\!<\!\sigma_{n-1}\!>\!\sigma_n$,}
\end{array}
\right.
\end{equation}
or, equivalently
\begin{equation}
\label{FSigreccb}
\F_\sigma(\X) = 
\partial\F_\sigma(\X) \cdot
\frac{q^{(n-1)a}-q^{(n-2)b}t}{1-q^{(n-2)b}t}\,\,
\frac{1-q^{(n-1)(1-a)}t}{1-q^n},
\end{equation}
where $a=1$  if $\sigma_{n-1}\!>\!\sigma_n$ and $a=0$ otherwise,
and  $b=1$  if $\sigma_{n-2}\!<\!\sigma_{n-1}\!>\!\sigma_n$
and $b=0$ otherwise.
\end{corollary}

As before, this result can be represented graphically with analogs of
the hook-content factors, by placing into node $i$ (from bottom to top) the 
$i$-th factor of $\F_\sigma(\X)$ of~(\ref{FSigrecc}).
For example, the second tree of Figure~\ref{figEqsFQS} shows the expansion of
$\F_\sigma(\X)$ with our second hook-contents of
$\sigma=(5, 6, 7, 4, 3, 2, 8, 9, 10, 1, 11)$.

The hook-content factors described in Corollary~\ref{propFsig} can have either
two or four terms.
But one easily checks that, if a factor has four terms, those 
terms simplify with the factors associated to the preceding letter in the
permutation. We  recover in this way the partial factors
of~\cite{NCSF2} and obtain a third version of the hook-content formula:

\begin{corollary}[\cite{NCSF2}, (152)]
\label{FSigNCSF2}
\begin{equation}
\label{FSigcorr}
\F_\sigma(\X) =
\prod_{i=1}^n
\frac{1}{1-q^i}
\left\{
\begin{array}{ll}
1-q^{i-1}t
& \text{\rm if $\sigma_{i-1}\!<\!\sigma_i\!<\!\sigma_{i+1}$,} \\[.1cm]
1-t
& \text{\rm if $\sigma_{i-1}\!<\!\sigma_i\!>\!\sigma_{i+1}$,} \\[.1cm]
q^{i-1}-t
& \text{\rm if $\sigma_{i-2}\!>\!\sigma_{i-1}\!>\!\sigma_i$}, \\[.1cm]
q^{i-1}-q^{i-2}t
& \text{\rm if $\sigma_{i-2}\!<\!\sigma_{i-1}\!>\!\sigma_i$,}
\end{array}
\right.
\end{equation}
with the conventions $\sigma_0=0$ and $\sigma_{n+1}=+\infty$.
\end{corollary}

The third tree of Figure~\ref{figEqsFQS} shows the 
resulting expansion of $\F_\sigma(\X)$  for the permutation
$\sigma=(5, 6, 7, 4, 3, 2, 8, 9, 10, 1, 11)$.
Note that it is obtained by permuting cyclically the numerators of the first 
formula among right branches.

%%%%%%%%%%%%%%%%%%%%%%%%%%%%%%%%%%%%%%%%%%%%%%%%%%%%%%%%%%%%%%%%%%%%%%%%%%%%%%%
%%%%%%%%%%%%%%%%%%%%%%%%%%%%%%%%%%%%%%%%%%%%%%%%%%%%%%%%%%%%%%%%%%%%%%%%%%%%%%%
%%%%%%%%%%%%%%%%%%%%%%%%%%%%%%%%%%%%%%%%%%%%%%%%%%%%%%%%%%%%%%%%%%%%%%%%%%%%%%%
\section{Compatibility between the dendriform operations and 
specialization of the alphabet}

%%%%%%%%%%%%%%%%%%%%%%%%%%%%%%%%%%%%%%%%%%%%%%%%%%%%%%%%%%%%%%%%%%%%%%%%%%%%%%%
\subsection{Dendriform algebras}

A dendriform algebra~\cite{Lod} is an associative algebra whose multiplication
$\cdot$  splits into two operations
\begin{equation}
a\cdot b = a \gaudend b + a \droitdend b
\end{equation}
satisfying
\begin{equation}
\left\{
\begin{array}{rcl}
(x\gaudend y)\gaudend z   &=& x\gaudend (y\cdot z)\,,\\
(x\droitdend y)\gaudend z &=& x\droitdend (y\gaudend z)\,,\\
(x\cdot y)\droitdend z    &=& x\droitdend (y\droitdend z)\,.
\end{array}
\right.
\end{equation}
For example, $\FQSym$ is dendriform with the following rules
\begin{equation}
\G_\alpha \gaudend \G_\beta
= \sum_{\gf{\gamma=uv\in \alpha* \beta}{|u|=|\alpha| ;\, \max(v)<\max(u)}}
  \G_\gamma\,,
\end{equation}
\begin{equation}
\G_\alpha \droitdend \G_\beta
= \sum_{\gf{\gamma=uv\in \alpha* \beta}{|u|=|\alpha| ;\, \max(v)\geq\max(u)}}
  \G_\gamma\,.
\end{equation}
Note that $x=\G_1=\F_1$ generates a free dendriform dialgebra in $\FQSym$,
isomorphic to $\PBT$, the Loday-Ronco algebra of planar binary
trees~\cite{LR1}.

%%%%%%%%%%%%%%%%%%%%%%%%%%%%%%%%%%%%%%%%%%%%%%%%%%%%%%%%%%%%%%%%%%%%%%%%%%%%%%%
\subsection{The half-products and the specialization}

\subsubsection{Descent statistics on half-shuffles}

On the basis $\F_\sigma$, the half-products are shifted half-shuffles.
Recall that the half-shuffles are the two terms of the recursive definition
of the shuffle product. For an alphabet $A$, and two words $u=u'a$,
$v=v'b$, $a,b\in A$, one has
\begin{equation}
u\shuffle v =u\gaudend v+u\droitdend v\,,
\end{equation}
where
\begin{equation}
u\gaudend v= (u'\shuffle v)a\ \text{and}\
u\droitdend v = (u\shuffle v')b\,.
\end{equation}
Assuming now that $A$ is totally ordered, we want to investigate
the distribution of descents on half-shuffles. To this aim we introduce
a linear map
\begin{equation}
\<w\>=F_{\D(w)}(X)=\<w|\sigma_1(XA)\>
\end{equation}
from $\K\<A\>$ to $QSym(X)$, the scalar product on $\K\<A\>$ being
defined by $\<u|v\>=\delta_{u,v}$.

For $w\in A^*$, let $\alph(w)\subseteq A$ be the set of letters occuring in
$w$.

\begin{lemma}
\label{desshuf}
If $\alph(u)\cap\alph(v)=\emptyset$, then 
\begin{equation}
\<u\shuffle v\>=\<u\>\<v\>\,.
\end{equation}
In particular, the descents of the elements of a shuffle on disjoint alphabets
depend only on the descents of the initial elements.
\end{lemma}

\Proof Denote by $\Delta$ the canonical (unshuffle) coproduct of $\K\<A\>$,
and write $u'v''$ for $u\otimes v$, so that $\Delta(a)=1\otimes a+a\otimes 1=a'+a''$
for $a\in A$. Then,
\begin{equation}
\begin{split}
\<u\shuffle v\>
&=\<u\shuffle v|\sigma_1(XA)\>
 =\left\<u' v''|\prod_{x\in X}^\rightarrow \Delta\sigma_x(A)\right\> \\
&=\left\<u' v''|\prod_{x\in X}^\rightarrow\prod_{a\in A}^\rightarrow
(1-x(a'+a''))^{-1}\right\> \\
&=\left\<u' v''|\prod_{x\in X}^\rightarrow\prod_{a'\in\alph(u')}^\rightarrow
(1-xa')^{-1}
\prod_{a''\in\alph(v'')}^\rightarrow (1-xa'')^{-1}\right\> \\
&=\<u|\sigma_1(XA)\>\<v|\sigma_1(XA)\>=\<u\>\<v\>\,.
\end{split}
\end{equation}
\qed

There is an equivalent statement for the dendriform half-products.

\begin{theorem}
\label{deshalf}
Let $u=u_1\cdots u_k$ and $v=v_1\cdots v_l$ of respective lengths $k$ and $l$.
If $\alph(u)\cap\alph(v)=\emptyset$, then
\begin{equation}
\<u\gaudend v\>=\<\sigma\gaudend \tau\>
\end{equation}
where $\sigma=\std(u)$ and $\tau=\std(v)[k]$ if $u_k<v_l$,
and $\sigma=\std(u)[l]$ and $\tau=\std(v)$ if $u_k>v_l$.
\end{theorem}

\Proof It is enough to check the first case, so we assume $u_k<v_l$.
The proof proceeds by induction on $n=k+l$. Let us set
$u=u'a'a$ and $\std(u)=u'_1a'_1a_1$.

If $a'>a'$, since $u\gaudend v=(u'a'\shuffle v)a$, we have
\begin{equation}
\<u\gaudend v\>=\sum_{w\in u'a'\shuffle v}F_{\D(w)\cdot 1}
=\<(u'_1a'_1\shuffle\tau)\cdot a_1\>
\end{equation}
with $\tau=\std(v)[k]$, according to Lemma \ref{desshuf}.

If $a'<a$, write $u\gaudend v=(u'a'\gaudend v)\cdot a+(u'a'\droitdend v)\cdot
a$. From the induction hypothesis, we have, with  $\tau$ as above,
$\<u'a'\gaudend v\>=\<u'_1a'_1\gaudend \tau\>$ and
$\<u'a'\droitdend v\>=\<u'_1a'_1\droitdend \tau\>$, so that
\begin{equation}
\<u\gaudend v\>=\sum_{w\in u'_1a'_1\gaudend\tau}F_{\D(w)\triangleright 1}
+\sum_{w\in u'_1a'_1\droitdend\tau}F_{\D(w)\cdot 1}\,,
\end{equation}
as required.
\qed

For example,
\begin{equation}
\begin{split}
\<634\gaudend 125\> =& \<631254 + 613254 + 612354 + 612534 + 163254 \\
+& 162354 + 162534 + 126354 + 126534 + 125634\>,
\end{split}
\end{equation}
\begin{equation}
\begin{split}
\<312\gaudend 456\> =& \<314562 + 341562 + 345162 + 345612 + 431562 \\
+& 435162 + 435612 + 453162 + 453612 + 456312\>,
\end{split}
\end{equation}
and one can check that both expressions are equal to
\begin{equation}
F_{132} + F_{141} + F_{1131} + F_{1221} + F_{222} + F_{231} + F_{2121}
+ F_{312} + F_{321} + F_{42}.
\end{equation}

\begin{corollary}
\label{deshalfCor}
Let $u$ and $v$ be two words of respective lengths $k$ and $l$.
Then, if $\alph(u)\cap\alph(v)=\emptyset$,
\begin{equation}
\label{egalMajs}
\sum_{x\in u\gaudend v} q^{\maj(x)} =
\sum_{y\in \sigma\gaudend \tau} q^{\maj(y)}.
\end{equation}
where $\sigma=\std(u)$ and $\tau=\std(v)[k]$ if $u_k<v_l$,
and $\sigma=\std(u)[l]$ and $\tau=\std(v)$ if $u_k>v_l$.
\end{corollary}

%%%%%%%%%%%%%%%%%%%%%%%%%%%%%%%%%%%%%%%%%%%%%%%%%%%%%%%%%%%%%%%%%%%%%%%%%%%%%%%
\subsubsection{$(q,t)$-specialization}
\label{qtDend}

We shall now see that Theorem~\ref{deshalf} implies a hook-content
formula for half-products evaluated over $\X$.
Let $\sigma\in\SG_n$ and $\tau\in\SG_m$.
Recall that $\tau[n]$ denotes the word $\tau_1+n,\tau_2+n,\dots,\tau+m+n$.
We have
\begin{equation}
\begin{split}
(q)_{n+m} \left(\F_\sigma \gaudend \F_\tau\right)(\X)
&=
\sum_{\epsilon\in\{\pm 1\}^{n+m}}
\sum_{\mu\in \sigma\gaudend \tau[n]} (-t)^{m(\epsilon)} q^{\maj(\mu,\epsilon)} \\
& = 
\sum_{\gf{\epsilon_1\in\{\pm 1\}^{n}}{\epsilon_2\in\{\pm 1\}^{m}}} 
\sum_{\mu'\in (\sigma,\epsilon_1)\gaudend (\tau[n],\epsilon_2)} (-t)^{m(\epsilon)} q^{\maj(\mu')},
\end{split}
\end{equation}
where $(\sigma,\epsilon_1)$, $(\tau[n],\epsilon_2)$, and $\mu'$ are signed
words.
Then, thanks to Theorem~\ref{deshalf}, the inner sum is the generating
function of the $\maj$ statistic on the left dendriform product of two
permutations. Its value is known (see~\cite{HNT2}, Equation~(34)), and is
\begin{equation}
q^{\maj(\sigma,\epsilon_1)} q^{\maj(\tau,\epsilon_2)} C(q),
\end{equation}
where $C(q)$ only depends on the sizes of $\sigma$ and $\tau$.

This implies that, if $\epsilon_1$ and $\epsilon'_1$ are two sign words
differing only on one entry,

\begin{equation}
\sum_{\mu'\in (\sigma,\epsilon_1)\gaudend (\tau[n],\epsilon_2)} q^{\maj(\mu')}
\text{\ and\ }
\sum_{\mu'\in (\sigma,\epsilon'_1)\gaudend (\tau[n],\epsilon_2)} q^{\maj(\mu')}
\end{equation}
are equal up to a power of $q$. Moreover, thanks to Lemma~\ref{lemFsig}, this
factor depends only on $\sigma$ and is the same as in the Lemma. The same
holds for two sign words $\epsilon_2$ and $\epsilon'_2$ differing on one entry, except 
for the last value of $\tau[n]$. In that special case, the contribution of the
letter is not given by Lemma~\ref{lemFsig} but by a similar statement where
the convention $y_{m}=0$ is replaced by $y_{m}=1$.
Hence, we have, deducing the second formula from the first one, since their
sum is $\F_\sigma(\X)\F_\tau(\X)$:

\begin{corollary}\label{cordend}
Let $\sigma\in\SG_n$ and $\tau\in\SG_m$. Then
\begin{equation}
\left(\F_\sigma \gaudend \F_\tau\right)(\X)
=
\frac{1-q^n}{1-q^{n+m}}\,\,
\frac{q^m-q^{(m-1)d}t}{1-q^{(m-1)d}t}\,\,
\F_{\sigma}(\X) \F_{\tau}(\X),
\end{equation}
and
\begin{equation}
\left(\F_\sigma \droitdend \F_\tau\right)(\X)
=
\frac{1-q^m}{1-q^{n+m}}\,\,
\frac{1-q^{n+(m-1)d}t}{1-q^{(m-1)d}t}\,\,
\F_{\sigma}(\X) \F_{\tau}(\X),
\end{equation}
where $d$ is  $1$ if $\tau_{m-1}<\tau_m$ and $0$ otherwise.
\end{corollary}

\begin{example}{\rm
Let us present all possible cases on the left dendriform product.
\begin{equation}
\begin{split}
(\F_{3421}\gaudend\F_{132})(\X)
&= q^2 \frac{(q-t)^2(1-t)^3(q^3-t)^2}{(1-q^7)(1-q^3)^2(1-q^2)^2(1-q)^2}\\
&= \frac{1-q^4}{1-q^7}\,\, \frac{q^3-t}{1-t} \,
\F_{3421}(\X) \F_{132}(\X).
\end{split}
\end{equation}
\begin{equation}
(\F_{3241}\gaudend\F_{213})(\X)
= \frac{1-q^4}{1-q^7}\,\,\frac{q^3-q^2t}{1-q^2t}\,
\F_{3241}(\X) \F_{213}(\X).
\end{equation}
\begin{equation}
(\F_{25134}\gaudend\F_{3421})(\X)
= \frac{1-q^5}{1-q^9}\,\,\frac{q^4-t}{1-t}\,
\F_{25134}(\X) \F_{3421}(\X).
\end{equation}
}
\end{example}

%%%%%%%%%%%%%%%%%%%%%%%%%%%%%%%%%%%%%%%%%%%%%%%%%%%%%%%%%%%%%%%%%%%%%%%%%%%%%%%
%%%%%%%%%%%%%%%%%%%%%%%%%%%%%%%%%%%%%%%%%%%%%%%%%%%%%%%%%%%%%%%%%%%%%%%%%%%%%%%
%%%%%%%%%%%%%%%%%%%%%%%%%%%%%%%%%%%%%%%%%%%%%%%%%%%%%%%%%%%%%%%%%%%%%%%%%%%%%%%
\section{A hook-content formula for binary trees}

%%%%%%%%%%%%%%%%%%%%%%%%%%%%%%%%%%%%%%%%%%%%%%%%%%%%%%%%%%%%%%%%%%%%%%%%%%%%%%%
\subsection{Subalgebras of $\FQSym$}

Recall that $\PBT$, the Loday-Ronco algebra of planar binary trees~\cite{LR1},
is naturally a subalgebra of $\FQSym$, the embedding being
\begin{equation}
\P_T(A)=\sum_{P(\sigma)=T}\F_\sigma(A)\,,
\end{equation}
where $P(\sigma)$ is the shape of the binary search tree associated with
$\sigma$~\cite{HNT}. Hence, $\P_T(\X)$ is well defined.

It was originally defined~\cite{LR1} as the free dendriform algebra on one
generator as follows: if $T$ is a binary tree $T_1$ (resp. $T_2$) be its left
(resp. right) subtree, then
\begin{equation}
\label{defPBT}
\P_T = \P_{T_1} \droitdend \P_1 \gaudend \P_{T_2}.
\end{equation}

%%%%%%%%%%%%%%%%%%%%%%%%%%%%%%%%%%%%%%%%%%%%%%%%%%%%%%%%%%%%%%%%%%%%%%%%%%%%%%%
\subsection{Hook-content formulas in PBT}

Note first that Corollary \ref{cordend} implies that the left and right
dendriform half-products factorize in the $\X$-specialization.
Because of the different expressions on signed permutations, it also proves
that the same property holds for trees, thanks to~(\ref{defPBT}).

Then, as a corollary of the definition of $\P_T$ and Corollary~\ref{corFAt},
we have

\begin{corollary}
Let $T$ be a binary tree. Then
\begin{equation}
\P_{T}(\X) = \frac{1}{(q)_n} \sum_{(\sigma,\epsilon)| P(\sigma)=T}
              (-t)^{m(\epsilon)} q^{\maj(\sigma,\epsilon)}.
\end{equation}
\end{corollary}

Recall that any binary tree has a unique standard labelling that makes it a
binary search tree. We then define the hook-content of a given node as the
contribution of its label among all signed permutations having this tree as
binary search tree. Thanks to Corollary \ref{cordend}, we get
a two-parameter version of the $q$-hook-length formulas
of Bj\"orner and Wachs \cite{BW1,BW2} (see also \cite{HNT2}):

\begin{theorem}
\label{PT1}
Let $T$ be a tree and $s$ a node of $T$. Let $n$ be the size of the subtree
whose root is $s$.
The $(q,t)$-hook-content factor of $s$ into $T$ is given by the following
rules:
\begin{equation}
\label{EqArb1}
\frac{1}{1-q^n}
\left\{
\begin{array}{lll}
q^{n} - q^{n'}t & \text{if $s$ is the right son of its father,}\\
1- q^{n'}t & \text{otherwise},
\end{array}
\right.
\end{equation}
where $n'$ is the size of the left subtree of $s$.
\end{theorem}

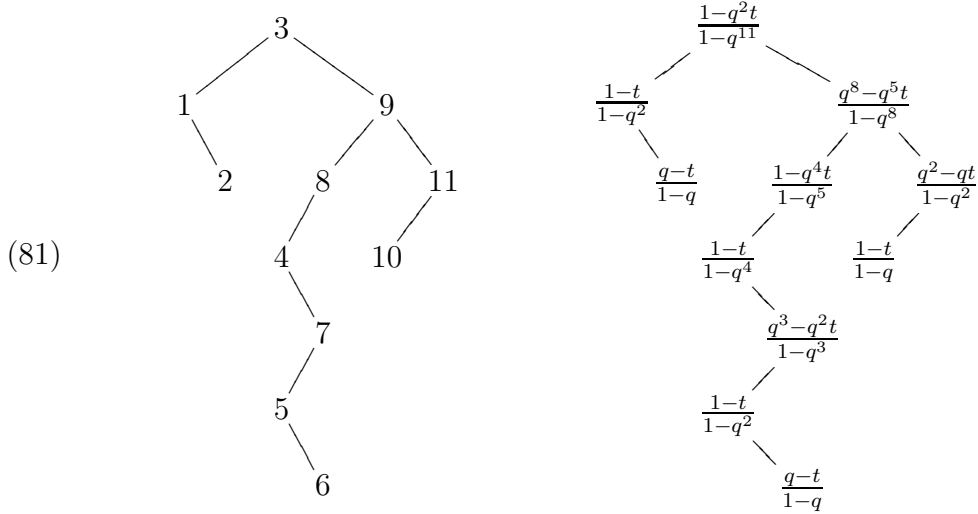
\begin{figure}[ht]
\entrymodifiers={+<4pt>}
\begin{equation}
\vcenter{\xymatrix@C=2mm@R=+6mm{
%etage 0
*{} & *{} & *{}
& 3 \ar@{-}[drrr]\ar@{-}[dlll] \\
%etage 1
  1 \ar@{-}[dr]
& *{} & *{} & *{} & *{} & *{} 
& 9 \ar@{-}[dll]\ar@{-}[dr] \\
%etage 2
*{} 
& 2  & *{} & *{}
& 8  \ar@{-}[dl] & *{} & *{}
& 11 \ar@{-}[dl]\\
%etage 3
*{} & *{} &*{}
& 4 \ar@{-}[dr]
& *{} & *{}
& 10 \\
%etage 4
*{} & *{} & *{} & *{}
& 7 \ar@{-}[dl]\\
%etage 5
*{} & *{} & *{}
& 5 \ar@{-}[dr]\\
%etage 6
*{} & *{} & *{} & *{}
& 6 \\
      }}
\qquad\qquad
\vcenter{\xymatrix@C=-1mm@R=+3mm{
%etage 0
*{} & *{} & *{}
& \frac{1-q^2t}{1-q^{11}}\ar@{-}[drrr]\ar@{-}[dlll] \\
%etage 1
\frac{1-t}{1-q^2}\ar@{-}[dr]
& *{} & *{} & *{} & *{} & *{} 
& \frac{q^8-q^5t}{1-q^8}\ar@{-}[dll]\ar@{-}[dr] \\
%etage 2
*{} 
& \frac{q-t}{1-q} & *{} & *{}
& \frac{1-q^4t}{1-q^5}\ar@{-}[dl] & *{} & *{}
& \frac{q^2-qt}{1-q^2}\ar@{-}[dl]\\
%etage 3
*{} & *{} &*{}
& \frac{1-t}{1-q^4}\ar@{-}[dr]
& *{} & *{}
& \frac{1-t}{1-q} \\
%etage 4
*{} & *{} & *{} & *{}
& \frac{q^3-q^2t}{1-q^3}\ar@{-}[dl]\\
%etage 5
*{} & *{} & *{}
& \frac{1-t}{1-q^2}\ar@{-}[dr]\\
%etage 6
*{} & *{} & *{} & *{}
& \frac{q-t}{1-q}\\
      }}
\end{equation}
\caption{\label{fig-PBT2}A binary tree (left diagram) labelled as a standard
binary search tree and the first hook-content formula on trees
(right diagram).
}
\end{figure}

As in the case of $\FQSym$, this  can be represented graphically by placing
into each node the fraction appearing in Equation~(\ref{EqArb1}).
For example, the first tree of Figure~\ref{figEqsFQS} shows the expansion of
$\F_\sigma(\X)$ with the first hook-contents of a zig-zag tree.
Figure~\ref{fig-PBT2} gives another example of this construction.

In particular, replacing $t$ by $-t$ in all formulas, this gives the following
combinatorial interpretation of the $(q,t)$ hook-length formula (recall that
$P(\sigma)={\mathcal T}(\sigma^{-1})$, where ${\mathcal T}(\tau)$ denotes the
decreasing tree of $\tau$):

\begin{corollary}
Let $T$ be a binary tree. Then the generating function of signed permutations
of shape $T$ by major index and number of signs is:
\begin{equation}
(q)_n \P_{T}(\X)|_{t=-t}
= \sum_{(\sigma,\epsilon)| P(\sigma)=T}
  t^{m(\epsilon)} q^{\maj(\sigma,\epsilon)}.
\end{equation}
\end{corollary}

\begin{example}
{\rm 
For example, with
\begin{equation}
T =
\vcenter{\xymatrix@C=2mm@R=+6mm{
%etage 0
*{} & *{} & *{}
& 3 \ar@{-}[drrr]\ar@{-}[dlll] \\
%etage 1
  1 \ar@{-}[dr]
& *{} & *{} & *{} & *{} & *{} 
& 5 \ar@{-}[dll]\ar@{-}[dr] \\
%etage 2
*{} 
& 2  & *{} & *{}
& 4  & *{} & *{}
& 6 \\
}},
\end{equation}
one has:
\begin{equation}
\begin{split}
\sum_{(\sigma,\epsilon)| P(\sigma)=T} t^{m(\epsilon)}q^{\maj(\sigma,\epsilon)}
& = (q)_6 \frac{(q+t)^2 (1+t)^2 (1+q^2t)
(q^3+qt)}{(1-q)^3(1-q^2)(1-q^3)(1-q^6)} \\
&= (q+t)^2 (1+t)^2 (1+q^2t) [4]_q [5]_q
\end{split}
\end{equation}
}
\end{example}

Here are the analogs of the other two hook-content formulas of~$\FQSym$.

\begin{theorem}
Let $T$ be a binary tree and $T_1$ (resp. $T_2$) be its left (resp. right)
subtree. Let $T'_2$ be the left subtree of $T_2$.
We then have
\begin{equation} 
\label{PTrecc}
\P_T(\X) =
\frac{(q^{\#T_2} - q^{\#T'_2}t)(1-q^{\#T_1})}{(1-q^{\#T'_2}t)(1-q^n)}
%\frac{1-q^{\#T_1}}{1-q^n}
\P_{T_1}(\X) \P_{T_2}(\X).
\end{equation} 
\end{theorem}

\Proof
This is a direct consequence of the dendriform specializations in $\FQSym$
thanks to~(\ref{defPBT}).
\qed

As in the case of $\FQSym$, it is possible to simplify the product 
so as to obtain a single quotient at each node.

\begin{corollary}
\label{CorPT}
Let $T$ be a tree and $s$ a node of $T$. Let $n$ be the size of the subtree
whose root is $s$.
The $(q,t)$-hook-content factor of $s$ into $T$ is given by the following rules:
\begin{equation}
\label{EqArb}
\frac{1}{1-q^n}
\left\{
\begin{array}{lll}
q^{n'} - q^{n''}t & \text{if $s$ has a right son,}\\
1- q^{n-1}t & \text{if $s$ has no right son and is not the right son of its
            father,}\\
1- q^{d}t & \text{if $s$ has no right son and is the right son of its
            father,}
\end{array}
\right.
\end{equation}
where $n'$ is the size of the right subtree of $s$,
$n''$ is the size of the left subtree of the right subtree of $s$,
and $d$ is the size where of the left subtree of the topmost ancestor of $s$
leading to $s$ only by right branches.
\end{corollary}

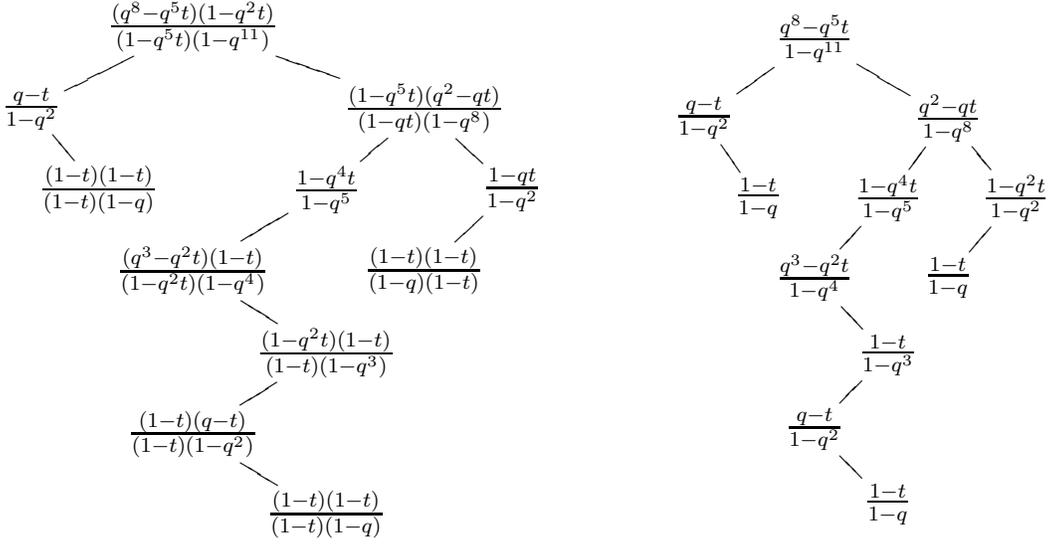
\begin{figure}[ht]
\entrymodifiers={+<4pt>}
\begin{equation*}
\vcenter{\xymatrix@C=-4mm@R=+3mm{
%etage 0
*{} & *{} & *{}
& \frac{(q^8-q^5t)(1-q^2t)}{(1-q^5t)(1-q^{11})}\ar@{-}[drrr]\ar@{-}[dlll] \\
%etage 1
\frac{q-t}{1-q^2}\ar@{-}[dr]
& *{} & *{} & *{} & *{} & *{} 
& \frac{(1-q^5t)(q^2-qt)}{(1-qt)(1-q^8)}\ar@{-}[dll]\ar@{-}[dr] \\
%etage 2
*{} 
& \frac{(1-t)(1-t)}{(1-t)(1-q)} & *{} & *{}
& \frac{1-q^4t}{1-q^5}\ar@{-}[dl] & *{} & *{}
& \frac{1-qt}{1-q^2}\ar@{-}[dl]\\
%etage 3
*{} & *{} &*{}
& \frac{(q^3-q^2t)(1-t)}{(1-q^2t)(1-q^4)}\ar@{-}[dr]
& *{} & *{}
& \frac{(1-t)(1-t)}{(1-q)(1-t)} \\
%etage 4
*{} & *{} & *{} & *{}
& \frac{(1-q^2t)(1-t)}{(1-t)(1-q^3)}\ar@{-}[dl]\\
%etage 5
*{} & *{} & *{}
& \frac{(1-t)(q-t)}{(1-t)(1-q^2)}\ar@{-}[dr]\\
%etage 6
*{} & *{} & *{} & *{}
& \frac{(1-t)(1-t)}{(1-t)(1-q)}\\
      }}
\qquad\qquad
%\begin{equation}
\vcenter{\xymatrix@C=-1mm@R=+3mm{
%etage 0
*{} & *{} & *{}
& \frac{q^8-q^5t}{1-q^{11}}\ar@{-}[drrr]\ar@{-}[dlll] \\
%etage 1
\frac{q-t}{1-q^2}\ar@{-}[dr]
& *{} & *{} & *{} & *{} & *{} 
& \frac{q^2-qt}{1-q^8}\ar@{-}[dll]\ar@{-}[dr] \\
%etage 2
*{} 
& \frac{1-t}{1-q} & *{} & *{}
& \frac{1-q^4t}{1-q^5}\ar@{-}[dl] & *{} & *{}
& \frac{1-q^2t}{1-q^2}\ar@{-}[dl]\\
%etage 3
*{} & *{} &*{}
& \frac{q^3-q^2t}{1-q^4}\ar@{-}[dr]
& *{} & *{}
& \frac{1-t}{1-q} \\
%etage 4
*{} & *{} & *{} & *{}
& \frac{1-t}{1-q^3}\ar@{-}[dl]\\
%etage 5
*{} & *{} & *{}
& \frac{q-t}{1-q^2}\ar@{-}[dr]\\
%etage 6
*{} & *{} & *{} & *{}
& \frac{1-t}{1-q}\\
      }}
\end{equation*}
\caption{\label{fig-PBT} Second and third hook-content formulas of a binary
tree:  by induction (left diagram),  simplification of the
induction (right diagram).}
\end{figure}

For example, on Figure~\ref{fig-PBT}, the rightmost node of the second tree
has coefficient $\frac{1-q^2t}{1-q^2}$: its topmost ancestor is the root
of the tree and the left subtree of the root is of size~$2$.
Note that it is obtained by permuting cyclically the numerators of the first 
formula among right branches, as it was already the case in $\FQSym$.

%%%%%%%%%%%%%%%%%%%%%%%%%%%%%%%%%%%%%%%%%%%%%%%%%%%%%%%%%%%%%%%%%%%%%%%%%%%%%%%
%%%%%%%%%%%%%%%%%%%%%%%%%%%%%%%%%%%%%%%%%%%%%%%%%%%%%%%%%%%%%%%%%%%%%%%%%%%%%%%
%%%%%%%%%%%%%%%%%%%%%%%%%%%%%%%%%%%%%%%%%%%%%%%%%%%%%%%%%%%%%%%%%%%%%%%%%%%%%%%
\section{Word Super-quasi-symmetric functions}

%%%%%%%%%%%%%%%%%%%%%%%%%%%%%%%%%%%%%%%%%%%%%%%%%%%%%%%%%%%%%%%%%%%%%%%%%%%%%%%
\subsection{Word quasi-symmetric functions}

Recall that a word $u$ over $\NN^*$ is 
\emph{packed} if the set of letters
appearing in $u$ is an interval of $\NN^*$ containing $1$.
Recall also that $\WQSym(A)$ is defined as the subalgebra of $\K\<A\>$ indexed
by \emph{packed words} and spanned by the elements
\begin{equation}
\M_u(A) := \sum_{\pack(w)=u} w,
\end{equation}
where $\pack(w)$ is the \emph{packed word} of $w$, that is, the word obtained
by replacing all occurrences of the $k$-th smallest letter of $w$ by $k$.
For example,
\begin{equation}
\pack(871883319) = 431442215.
\end{equation}

Let $\NW_u=\M_u^*$ be the dual basis of $(\M_u)$.
It is known that $\WQSym$ is a self-dual Hopf algebra~\cite{Hiv,NT06} and that
on the dual $\WQSym^*$, an internal product $*$ may be defined by
\begin{equation}
\label{intWQ}
\NW_u * \NW_v = \NW_{\pack(u,v)},
\end{equation}
where the packing of biwords is defined with respect to the lexicographic
order, so that, for example,
\begin{equation}
\pack\left(\gf{42412253}{53154323}\right)
= 62513274.
\end{equation}

This product is induced from the internal product of parking
functions~\cite{NTpark,NT1,NTp2} and allows one to identify the homogeneous
components $\WQSym_n$ with the (opposite) Solomon-Tits algebras, in the sense
of~\cite{Patras}.

The (opposite) Solomon descent algebra, realized
as $\Sym_n$, is embedded in the (opposite) Solomon-Tits algebra realized as
$\WQSym^*_n$ by
\begin{equation}
S^I = \sum_{ev(u)=I} \NW_u.
\end{equation}

%%%%%%%%%%%%%%%%%%%%%%%%%%%%%%%%%%%%%%%%%%%%%%%%%%%%%%%%%%%%%%%%%%%%%%%%%%%%%%%
\subsection{An algebra on signed packed words}

Let us define $\WQSym^{(2)}$ as the space spanned by the
$\M_{u,\epsilon}$, where
\begin{equation}
\M_{u,\epsilon}(A) :=
 \sum_{\gf{(w,\epsilon)}{\pack(w)=u}} (w,\epsilon).
\end{equation}
This is  a Hopf algebra for the standard operations.
We denote by
$\NW_{u,\epsilon}$ the dual basis of $\M$.
This algebra contains $\NCSF^{(2)}$, the Mantaci-Reutenauer algebra of type
$B$. To show this, let us describe the embedding.

A signed word is said to be \emph{regular} if all occurences of any unsigned
letter have same sign.
For example, $11\ol2\ol231$ is regular, but $11\ol1$ and $1121\ol2$ are not.

The \emph{signed evaluation} $\sev(w,\epsilon)$ of a regular word
is the signed composition $(I,\mu)$ where $i_j$ is the number of occurrences
of the (unsigned) letter $j$ and $\mu_j$ is the sign of $j$ in $(w,\epsilon)$.

Let $\phi$ be the morphism from $\NCSF^{(2)}$ into $\WQSym^{(2)}$ defined by
\begin{equation}
\phi(S_n) = \NW_{1^n},
\qquad
\phi(S_{\ol n}) = \NW_{\ol1^n}.
\end{equation}

We then have :
\begin{lemma}
\begin{equation}
\phi(S^{(I,\epsilon)}) =
\sum_{\gf{(u,\epsilon') regular}{\sev(u,\epsilon')=(I,\epsilon)}}
\NW_{u,\epsilon'}.
\end{equation}
\end{lemma}

\Proof
This follows from the product formula of the $\NW$, which is a special
case of the multiplication of signed parking functions \cite{FQSl}.
\qed

The image of $\NCSF^{(2)}$ by this embedding is contained in the Hopf
subalgebra $\BW$ of $\WQSym^{(2)}$ generated by the $\NW_{u,\epsilon}$ indexed
by regular signed packed words.
The dimensions of the homogeneous components $\BW_n$ are given by
Sequence~A004123 of~\cite{Slo} whose first values are
\begin{equation}
1,\ 2,\ 10,\ 74,\ 730,\ 9002,\ 133210.
\end{equation}

Note in particular that $\sigma_1^\#$ has a simple expression in terms of
$\NW_{u,\epsilon}$.

\begin{lemma}
\label{sigmaOnN}
Let $PW$ denote the set of packed words, and $\max(u)$ the maximal letter of
$u$. Then
\begin{equation}
{\sigma_1^\#}
 = \sum_{u\in PW}
        (-1)^{n-\max(u)} \NW_{u,(-1)^n}
      + (-1)^{m(\epsilon')-(\max(u)-1)} \NW_{u,\epsilon'},
\end{equation}
where $(u,\epsilon')$ is such that all letters but the maximal one are
signed.
\end{lemma}

\begin{example}
\begin{equation}
S_2^\#=
  - \NW_{\ol1\ol1} + \NW_{11}
  + \NW_{\ol1\ol2} + \NW_{\ol12}
  + \NW_{\ol2\ol1} + \NW_{2\ol1}.
\end{equation}
\begin{equation}
\begin{split}
S_3^\# =&
+ \NW_{\ol1\ol1\ol1} + \NW_{111}\\
&
- \NW_{\ol1\ol1\ol2} - \NW_{\ol1\ol12}
- \NW_{\ol1\ol2\ol1} - \NW_{\ol12\ol1} - \NW_{\ol2\ol1\ol1}
- \NW_{2\ol1\ol1}\\
&
- \NW_{\ol2\ol2\ol1} + \NW_{   22\ol1} - \NW_{\ol2\ol1\ol2} + \NW_{2\ol1  2 }
- \NW_{\ol1\ol2\ol2} + \NW_{\ol12   2}\\
&+ \NW_{\ol1\ol2\ol3} + \NW_{\ol1\ol23}
+ \NW_{\ol1\ol3\ol2} + \NW_{\ol13\ol2}
+ \NW_{\ol2\ol1\ol3} + \NW_{\ol2\ol13}\\
& + \NW_{\ol2\ol3\ol1} + \NW_{\ol23\ol1}
+ \NW_{\ol3\ol1\ol2} + \NW_{3\ol1\ol2} + \NW_{\ol3\ol2\ol1} + \NW_{3\ol2\ol1}
.
\end{split}
\end{equation}
\end{example}

%%%%%%%%%%%%%%%%%%%%%%%%%%%%%%%%%%%%%%%%%%%%%%%%%%%%%%%%%%%%%%%%%%%%%%%%%%%%%%%
\subsection{An internal product on signed packed words}

The internal product of $\WQSym^*$~(\ref{intWQ}) can be extended to
${\WQSym^{(2)}}^*$ by
\begin{equation}
\label{intWQc}
\NW_{u,\epsilon} * \NW_{v,\rho} = \NW_{\pack(u,v),\epsilon\rho},
\end{equation}
where $\epsilon\rho$ is the componentwise product.
One obtains in this way the (opposite) Solomon-Tits algebra of type $B$.
This product is induced from the internal product of signed parking
functions \cite{FQSl} and can be shown to coincide with the one introduced
by Hsiao \cite{Hsiao}.

From this definition, we immediately have:
\begin{proposition}
$\BW$ is a subalgebra of ${\WQSym^{(2)}}^*$ for the internal product.
\end{proposition}

Since $\sigma_1^\#$ belongs to ${\WQSym^{(2)}}^*$, we can define
\begin{equation}
\NW_{u}^\# := \NW_{u}(A|\Ab) = \NW_u * \sigma_1^\#.
\end{equation}

\begin{example}
\label{exWQ}
Let us compute the first $\NW_u(A|\Ab)$.
\begin{equation}
\NW_{11}^\# = 
  - \NW_{\ol1\ol1} + \NW_{11}
  + \NW_{\ol1\ol2} + \NW_{\ol12}
  + \NW_{\ol2\ol1} + \NW_{2\ol1}.
\end{equation}
\begin{equation}
\NW_{12}^\# = \NW_{\ol1\ol2} + \NW_{\ol12} + \NW_{1\ol2} + \NW_{12}.
\end{equation}
\begin{equation}
\NW_{21}^\# = \NW_{\ol2\ol1} + \NW_{\ol21} + \NW_{2\ol1} + \NW_{12}.
\end{equation}
\begin{equation}
\begin{split}
\NW_{112}^\# =&
  - \NW_{\ol1\ol1\ol2} - \NW_{\ol1\ol12} + \NW_{11\ol2} + \NW_{112}\\
& + \NW_{\ol1\ol2\ol3} + \NW_{\ol1\ol23} + \NW_{\ol12\ol3} + \NW_{\ol123}\\
& + \NW_{\ol2\ol1\ol3} + \NW_{\ol2\ol13} + \NW_{2\ol1\ol3} + \NW_{2\ol13}\\
\end{split}
\end{equation}
\begin{equation}
\begin{split}
\NW_{121}^\# =&
  - \NW_{\ol1\ol2\ol1} - \NW_{\ol12\ol1} + \NW_{1\ol21} + \NW_{121}\\
& + \NW_{\ol1\ol3\ol2} + \NW_{\ol13\ol2} + \NW_{\ol1\ol32} + \NW_{\ol132}\\
& + \NW_{\ol2\ol3\ol1} + \NW_{\ol23\ol1} + \NW_{2\ol3\ol1} + \NW_{23\ol1}\\
\end{split}
\end{equation}
\end{example}

In the light of the previous examples, let us say that a packed word $v$ is
\emph{finer} than a packed $u$, and write $v\finerW u$ if
$u$ can be obtained from $v$ by application of a nondecreasing map from
$\NN^*$ to $\NN^*$. Note that this definition is easy to describe on set
compositions: $u$ is then obtained by gluing together consecutive
parts of $v$.
For example, the words finer than $121$ are $121$, $132$, and $231$.

\begin{theorem}
\label{thm-WQAAb}
Let $u$ be a packed word.
Then
\begin{equation}
\NW_{u}^\# = \sum_{v\finerW u}
             \sum_{\epsilon}
              (-1)^{m(\epsilon)+m'(v,\epsilon) } \NW_{v,\epsilon}
\end{equation}
where $m'(v,\epsilon)$ is equal to the number of \emph{different} signed
letters of $v$ and
where the sum on $\epsilon$ is such that the words $(v,\epsilon)$ are regular
and such  that if more than two letters of $v$ go to the same letter of
$u$, all letters but the greatest are signed (the greatest can be
either signed or not).
In particular, the number of such $\epsilon$ for a given $v$
is equal to $2^{\max(u)}$, so is independent ov $v$.
\end{theorem}

\Proof
From the definitions of $\sigma_1^\#$ and of the packing algorithm, it is
clear that the words appearing on the expansion of $\NW_{u}^\#$ are exactly
the words given in the previous statement.

Moreover, the coefficient of a signed word $(w,\epsilon)$ in
$\sigma_1^\#$ is equal to the coefficient of any of its rearrangements (where
the signs follow their letter). Now, given a permutation $\sigma$ and two
words $u$ and $u'$ having a word $v$ as packed word, the packed word of
$u\cdot\sigma$ and $u'\cdot\sigma$ is $v\cdot\sigma$. So we can restrict
ourselves to compute $\NW_u^\#$ for all nondecreasing words $u$ since all the
other ones are obtained by permutation of the entries.

Assume now that $u$ is a nondecreasing word, and let us show that the
coefficient of $(v,\epsilon)$ in $\NW_{u}^\#$ is either $1$ or $-1$.
The only terms $\NW$ in $\sigma_1^\#$ that can yield $(v,\epsilon)$ when
multiplied on the left by $\NW_u$ are the signed words with negative entries
exactly as in $\epsilon$. Let $T_\epsilon$ denote this set.
Thanks to Lemma~\ref{sigmaOnN}, the $\NW$ appearing in the expansion of
${\sigma_1^\#}_n$ with negative signs at $k$ given slots are the following
packed words: all the elements of $PW_k$ at the negative slots and one letter
greater than all the others at the remaining slots.
In particular, the cardinality of $T_\epsilon$  depends only on $k$ and is
equal to $|PW_k|$. Since there is only one positive value for each element,
two words $w$ and $w'$ of $T_\epsilon$ give the same result by packing $(u,w)$
and $(u,w')$ if they coincide on the negative slots.

This  means that we can restrict ourselves to the special case where
$\epsilon=(-1)^n$ since the positive slot do not change the way of regrouping
the elements of $T_\epsilon$ to obtain $(v,\epsilon)$.
Now, the sign has been disposed of and we can concentrate on the packing
algorithm. The previous discussion shows that we only need to prove that,
given a word $v$ finer than a word $u$, the set $T$ of packed words $w$ such
that $\pack(u,w)=v$ satisfies the following property: if $t_d$ is the number
of elements of $T$ with maximum $d$, then
\begin{equation}
\sum_d (-1)^d t_d = \pm1.
\end{equation}

Thanks to the packing algorithm, we see that $T$ is the set of packed words
with (in)equalities coming from the values of $v$ at the places where $u$ have
equal letters. So $T$ is a set of packed words with (in)equalities between
adjacent places with no other relations. Hence, if $u$ has $l$ different
letters, $T$ is obtained as the product of $l$ quasi-monomial functions
$\M_w$. The conclusion of the proof comes from the following lemma. 
\qed

\begin{lemma}
Let $w_1,\dots,w_k$ be $k$ packed words with respective maximum letters
$a_1,\dots,a_k$.
Let $T$ be the set of packed words appearing in the expansion of
\begin{equation}
\M_{w_1} \dots \M_{w_k}.
\end{equation}
Then, if $t_d$ is the number of elements of $T$ with maximum $d$, then
\begin{equation}
\sum_d (-1)^d t_d = (-1)^{a_1+\dots+a_k}.
\end{equation}
\end{lemma}

\Proof
We only need to prove the result for $k=2$ since the other cases follow by
induction: compute $\M_{w_1}\dots\M_{w_{k-1}}$ and multiply this by
$\M_{w_k}$ to get the result.

Let us compute $\M_{w_1} \M_{w_2}$.
The number of words with maximum $a_1+a_2-d$ in this product is equal to
\begin{equation}
\label{td}
\binom{a_1}{d} \binom{a_1+a_2-d}{a_1}.
\end{equation}
Indeed, a word in $\M_{w_1} \M_{w_2}$ with maximum $a_1+a_2-d$ is completely
characterized by the
$d$ integers between $1$ and $a_1+a_2-d$ common to the prefix of size $|w_1|$
and the suffix of size $|w_2|$ of $w$,
by the $(a_1-d)$ integers only appearing in the prefix,
and the $(a_2-d)$ integers only appearing in the suffix, which hence gives the
enumeration formula
\begin{equation}
t_{a_1+a_2-d} = \frac{(a_1+a_2-d)!}{d! (a_1-d)! (a_2-d)!},
\end{equation}
equivalent to the previous one.

It remains to compute
\begin{equation}
\sum_{d} (-1)^{a_1+a_2-d} \binom{a_1}{d} \binom{a_1+a_2-d}{a_1},
\end{equation}
which is, with the usual notation for elementary and complete homogeneous
symmetric functions, understood as operators of the $\lambda$-ring $\ZZ$,
\begin{equation}
\begin{split}
& (-1)^{a_1+a_2} \sum_d (-1)^d e_{d}(a_1) h_{a_2-d}(a_1+1) \\
&= (-1)^{a_1+a_2} \sum_d h_d(-a_1) h_{a_2-d}(a_1+1) \\
&= (-1)^{a_1+a_2} h_{a_2}(-a_1+a_1+1)\\
&= (-1)^{a_1+a_2} h_{a_2}(1) = (-1)^{a_1+a_2}.
\end{split}
\end{equation}
\qed

This combinatorial interpretation of (\ref{td}) gives back in particular one
interpretation of the Delannoy numbers (sequence A001850 of \cite{Slo}) and
of their usual refinement (sequence A008288 of \cite{Slo}).

%%%%%%%%%%%%%%%%%%%%%%%%%%%%%%%%%%%%%%%%%%%%%%%%%%%%%%%%%%%%%%%%%%%%%%%%%%%%%%%
\subsection{Specializations}

The internal product of $\WQSym^*$ allows in particular to define
\begin{equation}
\NW_{u}((1-t)A) := \NW_{u}(A) * \sigma_1((1-t)A) = \eta_t(\NW_u),
\end{equation}
so that we have
\begin{equation}
\S^I ((1-t)A) = \sum_{\ev(u)=I} \NW_u((1-t)A).
\end{equation}

\begin{example}{\rm
Taking the same five examples as in Example~\ref{exWQ}, we get
\begin{equation}
\NW_{11}((1-t)A) = 
  (1-t^2) \NW_{11} - t(1-t) \NW_{12} - t(1-t) \NW_{21}.
\end{equation}
\begin{equation}
\NW_{12}((1-t)A) = (1-t)^2 \NW_{12}
\qquad
\text{and}
\qquad
\NW_{21}((1-t)A) = (1-t)^2 \NW_{21}.
\end{equation}
\begin{equation}
\NW_{112}((1-t)A)
 = (1-t)(1-t^2) \NW_{112} -t(1-t)^2 \NW_{123} -t(1-t)^2 \NW_{213}.
\end{equation}
\begin{equation}
\NW_{121}((1-t)A)
 = (1-t)(1-t^2) \NW_{121} -t(1-t)^2 \NW_{132} -t(1-t)^2 \NW_{231}.
\end{equation}
}
\end{example}

\begin{theorem}
\label{1mtWQ}
Let $u$ be a packed word.
Then
\begin{equation}
\NW_{u}((1-t)A) = \sum_{v\finerW u} 
              (-1)^{\max(v)-\max(u)} t^{f(u,v)}
              \prod_{k=1}^{\max(u)} (1-t^{g(u,v,k)})
              \ \ \NW_v(A).
\end{equation}
where, if one writes
\begin{equation}
\ev(u)=(i_1,\dots,i_p) \text{ and }
\ev(v)= ( (i_1^{(1)},\dots,i_1^{(q_1)}),\dots, (i_p^{(1)},\dots,i_p^{(q_p)})),
\end{equation}
then
\begin{equation}
f(u,v) := \sum_{k=1}^p \sum_{j=1}^{q_k-1} i_k^{(j)}
\qquad
\text{and}
\qquad
g(u,v,k) := i_k^{(q_k)}.
\end{equation}
\end{theorem}

\Proof
This is a direct consequence of Theorem~\ref{thm-WQAAb}.
\qed

%%%%%%%%%%%%%%%%%%%%%%%%%%%%%%%%%%%%%%%%%%%%%%%%%%%%%%%%%%%%%%%%%%%%%%%%%%%%%%%
\subsection{Duality}

By duality, one defines
\begin{equation}
\M_{u}(A\cdot (1-t)) := \eta_t^*(\M_u(A)),
\end{equation}
since
\begin{equation}
\sum_{u} \M_u(A\cdot (1-t))\otimes\NW_u(B) =
\sum_{u} \M_u(A)\otimes\NW_u((1-t)B).
\end{equation}

\begin{example}
\begin{equation}
\M_{11}(A\cdot (1-t)) = (1-t^2) \M_{11}(A).
\end{equation}
\begin{equation}
\M_{12}(A\cdot (1-t)) = - t(1-t) \M_{11}(A) + (1-t)^2 \M_{12}(A).
\end{equation}
\begin{equation}
\M_{21}(A\cdot (1-t)) = - t(1-t) \M_{11}(A) + (1-t)^2 \M_{21}(A).
\end{equation}
\begin{equation}
\M_{112}(A\cdot (1-t)) = (1-t)(1-t^2) \M_{112}(A) - t^2(1-t) \M_{111}(A).
\end{equation}
\begin{equation}
\M_{121}(A\cdot (1-t)) = (1-t)(1-t^2) \M_{121}(A) - t^2(1-t) \M_{111}(A).
\end{equation}
\begin{equation}
\begin{split}
\M_{123}(A\cdot (1-t)) =& (1-t)^3 \M_{123}(A) - t(1-t)^2\M_{112}(A) \\
                 &- t(1-t)^2 \M_{122}(A) + t^2(1-t)\M_{111}(A).
\end{split}
\end{equation}
\end{example}

Since the transition matrix from $\M(A\cdot (1-t))$ to $\M(A)$ is the
transpose of the transition matrix from $\NW((1-t)A)$ to $\NW(A)$,
we can obtain a simple combinatorial interpretation of $\M(A\cdot (1-t))$.

First, let us define the \emph{super-packed word} $v:=\spack(u,\epsilon)$
associated with a regular signed word $(u,\epsilon)$.
Let $f_\epsilon$ be the nondecreasing function sending $1$ to $1$ and each value
$i$ either to $f_\epsilon(i-1)$ if the value $i-1$ is signed in $\epsilon$ or to
$1+f_\epsilon(i-1)$ if not. Extend $f_\epsilon$ to a morphism of $A^*$.
Then $v=f_\epsilon(u)$.

For example,
\begin{equation}
\spack(5\ol1\ol2\ol135\ol4\ol46\ol1) = 2111122231.
\end{equation}

Let $[v,u]$ be the interval for the refinement order on words, that is,
the set of packed words $w$ such that $u\finerW w\finerW v$.  
\begin{proposition}
\label{Mu-col}
Let $u$ be a word.
Then
\begin{equation}
\M_u(A\cdot (1-t)) = \sum_{(u,\epsilon) \text{regular}}
                     (-1)^{m'(u,\epsilon)} t^{m(\epsilon)}
                    \sum_{w\in [\spack(u,\epsilon),u]} \M_{w}(A).
\end{equation}
\end{proposition}

\Proof
Observe that if a signed word $(u,\epsilon)$ appears in $\NW_w^\#$ then it also
appears in $\NW_v^\#$ for all $v\in [u,w]$.
The rest comes directly from Theorem~\ref{thm-WQAAb} and from the fact that
$\NW_{(u,\epsilon)}$ is sent to $(-t)^{m(\epsilon)} \NW_{u}$ when sending
$\ol A$ to $-tA$.
\qed

\begin{example}{\rm
\begin{equation}
\M_{21}(A.(1-t)) = (-t+t^2) (\M_{11}+\M_{21}) + (1-t) \M_{21}.
\end{equation}
\begin{equation}
\M_{112}(A.(1-t)) = (-t^2+t^3) (\M_{111}+\M_{112}) +(1-t) \M_{112}.
\end{equation}
\begin{equation}
\begin{split}
\M_{123}(A.(1-t)) =& (t^2-t^3) (\M_{111}+\M_{112}+\M_{122}+\M_{123}) \\
                   & +(-t+t^2) (\M_{112}+\M_{123})\\
                   & +(-t+t^2) (\M_{122}+\M_{123})\\
                   & +(1-t)     \M_{123}.
\end{split}
\end{equation}
}
\end{example}

When $A$ is a commutative alphabet $X$, this specializes to $M_I(X(1-t))$
where $I=\ev(u)$ and in particular, for $X=\frac{1}{1-q}$, we recover a result
of~\cite{NCSF2}:

\begin{theorem}[\cite{NCSF2}]
Let $u$ be a packed word of size $n$.
\begin{equation}
\M_u(\X) = M_I (\X)
= \frac{1-t^{i_p}}{1-q^n}
  \prod_{k=1}^{p-1} \frac{(q^{i_1+\dots+i_k}-t^{i_k})}{1-q^{i_1+\dots+i_k}}.
\end{equation}
where the composition $I=(i_1,\dots,i_p)$ is the evaluation of $u$. 
\end{theorem}

\Proof
From Proposition~\ref{Mu-col} giving a combinatorial interpretation of
$\M_u(A\cdot (1-t))$, we have:
\begin{equation}
\M_u(\X) = 
\sum_{(u,\epsilon) \text{regular}}
                     (-1)^{m'(u,\epsilon)} t^{m(\epsilon)}
                    \sum_{w\in [\spack(u,\epsilon),u]} \M_{w}(1/(1-q)).
\end{equation}
We now have to evaluate the sum of $\M_{w}(1/(1-q))$ over an interval of the
composition lattice. Thanks to Lemma~\ref{lem-1mq}  below, it is
equal to
\begin{equation}
\frac{q^{\maj(I)}}{(1-q^{k_1})(1-q^{k_1+k_2})\cdots(1-q^{k_1+k_2+\cdots+k_s})},
\end{equation}
where $I=\ev(\spack(u,\epsilon))$ and $K=\ev(u)$,
which implies the result.  \qed

\begin{lemma}
\label{lem-1mq}
Let $I$ and $K$ be two compositions of $n$ such that $K\finer I$. Then
\begin{equation}
\sum_{J\in [I,K]} M_J(1/(1-q)) =
\frac1{1-q^n}  \frac{q^{\maj(I)}}{\prod_{d\in\Des(K)} 1-q^d}.
\end{equation}
\end{lemma}

\Proof
We have
\begin{equation}
M_J(1/(1-q)) = \frac{1}{1-q^n} \prod_{d\in Des(J)} \frac{q^d}{1-q^d}.
\end{equation}
Factorizing by the common denominator of all these elements and by
$q^{\maj(K)}$, we have to evaluate
\begin{equation}
\sum_{D\subseteq \Des(K)\backslash \Des(I)}\prod_{d\in D}(1-q^d) q^{-d}
\end{equation}
which is equal to
\begin{equation}
\prod_{d\in \Des(K)/\Des(I)}(1-1+q^{-d})
= q^{-(\maj(K)-\maj(I))}.
\end{equation}
\qed

Putting together Proposition~\ref{Mu-col} and Lemma~\ref{lem-1mq}, one
obtains:
\begin{corollary}
\label{cor-1mq}
Let $u$ be a word of size $n$.
Then
\begin{equation}
(q)_{\ev(u)} \M_u(\X) = \sum_{(u,\epsilon) \text{regular}}
                       (-1)^{m'(u,\epsilon)} t^{m(\epsilon)}
                       q^{\maj(\spack(u,\epsilon))},
\end{equation}
where $(q)_I$ is defined as $(1-q^n)\prod_{d\in\Des(I)}(1-q^d)$.
\end{corollary}

\begin{corollary}
\label{cor-M1mq}
Let $u$ be a word of size $n$.
Then the generating function of signed permutations
of unsigned part $u$ by major index of their super-packed word and number of
signs is:
\begin{equation}
\sum_{(u,\epsilon) \text{regular}}
    t^{m(\epsilon)} q^{\maj(\spack(u,\epsilon))}
= (1+t^{i_p}) \prod_{k=1}^{p-1} (q^{i_1+\dots+i_k}+t^{i_k}).
\end{equation}
\end{corollary}

\begin{example}
{\rm
For example, with $u = (1,1,2,3,3,3,3,4,4)$, one has:
\begin{equation}
\begin{split}
\sum_{(u,\epsilon) regular}
t^{m(\epsilon)} q^{\maj(\spack(u,\epsilon))}
= (1+t^2) (q^2+t^2) (q^3+t) (q^7+t^4).
\end{split}
\end{equation}
}
\end{example}

%%%%%%%%%%%%%%%%%%%%%%%%%%%%%%%%%%%%%%%%%%%%%%%%%%%%%%%%%%%%%%%%%%%%%%%%%%%%%%%
%%%%%%%%%%%%%%%%%%%%%%%%%%%%%%%%%%%%%%%%%%%%%%%%%%%%%%%%%%%%%%%%%%%%%%%%%%%%%%%
%%%%%%%%%%%%%%%%%%%%%%%%%%%%%%%%%%%%%%%%%%%%%%%%%%%%%%%%%%%%%%%%%%%%%%%%%%%%%%%
\section{Tridendriform operations and the specialization of alphabet}

%%%%%%%%%%%%%%%%%%%%%%%%%%%%%%%%%%%%%%%%%%%%%%%%%%%%%%%%%%%%%%%%%%%%%%%%%%%%%%%
\subsection{Tridendriform structure of $\WQSym$}

A \emph{dendriform trialgebra} \cite{LRtri} is an associative algebra whose
multiplication
$\cdot$ splits into three pieces
\begin{equation}
x\cdot y = x\gautrid y + x\miltrid y + x\droittrid y\,,
\end{equation}
where $\miltrid$ is associative, and

\begin{equation}
(x\gautrid y)\gautrid z = x\gautrid (y\cdot z)\,,\ \
(x\droittrid y)\gautrid z = x\droittrid (y\gautrid z)\,,\ \
(x\cdot y)\droittrid z = x\droittrid (y\droittrid z)\,,\ \
\end{equation}
\begin{equation}
(x\droittrid y)\miltrid z = x\droittrid (y\miltrid z)\,,\ \
(x\gautrid y)\miltrid z = x\miltrid (y\droittrid z)\,,\ \
(x\miltrid y)\gautrid z = x\miltrid (y\gautrid z)\,.
\end{equation}

It has been shown in~\cite{NT-cras} that the augmentation ideal
$\K\langle A\rangle^+$ has a natural structure of dendriform trialgebra:
for two non empty words $u,v\in A^*$, we set
\begin{eqnarray}
u\gautrid v=\begin{cases} uv &
\text{if $\max(u)>\max(v)$}\cr 0 &\mbox{otherwise,} \end{cases}\\
u\miltrid v=\begin{cases} uv &
\text{if $\max(u)=\max(v)$}\cr 0 &\mbox{otherwise,} \end{cases}\\
u\droittrid v=\begin{cases} uv &
\text{if $\max(u)<\max(v)$}\cr 0 &\mbox{otherwise.} \end{cases}
\end{eqnarray}

$\WQSym^+$ is a sub-dendriform trialgebra of $\K\<A\>^+$,
the partial products being given by
\begin{equation}
\M_{w'} \gautrid \M_{w''} =
\sum_{w=u\cdot v\in w'*_W w'', |u|=|w'| ; \max(v)<\max(u)}
\M_w,
\end{equation}
\begin{equation}
\M_{w'} \miltrid \M_{w''} =
\sum_{w=u\cdot v\in w'*_W w'', |u|=|w'| ; \max(v)=\max(u)}
\M_w,
\end{equation}
\begin{equation}
\M_{w'} \droittrid \M_{w''} =
\sum_{w=u\cdot v\in w'*_W w'', |u|=|w'| ; \max(v)>\max(u)}
\M_w,
\end{equation}
where the \emph{convolution} $u' *_W u''$ of two packed words
is defined as
\begin{equation} 
u' *_W u'' = \sum_{v,w ;
u=v\cdot w\,\in\,\PW, \pack(v)=u', \pack(w)=u''} u\,.
\end{equation}

%%%%%%%%%%%%%%%%%%%%%%%%%%%%%%%%%%%%%%%%%%%%%%%%%%%%%%%%%%%%%%%%%%%%%%%%%%%%%%%
\subsection{Specialization of the partial products}

If $w$ is a packed word, let $\Nmax(w)$ be the number of maximal letters
of $w$ in $w$.

\begin{theorem}
Let $u_1\in\PW(n)$ and $u_2\in\PW(m)$.
Then
\begin{equation}
(\M_{u_1} \gautrid \M_{u_2})(\X)
=
\frac{1-q^n}{1-q^{n+m}}\,\,
\frac{q^m-t^{\Nmax(u_2)}}{1-t^{\Nmax(u_2)}}\,\,
\M_{u_1}(\X) \M_{u_2}(\X),
\end{equation}
\begin{equation}
(\M_{u_1} \miltrid \M_{u_2})(\X)
=
\frac{(1-q^n)(1-q^m)}{1-q^{n+m}}\,\,
\frac{1-t^{\Nmax(u_1)+\Nmax(u_2)}}{(1-t^{\Nmax(u_1)})(1-t^{\Nmax(u_2)})}
\M_{u_1}(\X) \M_{u_2}(\X),
\end{equation}
and
\begin{equation}
(\M_{u_1} \droittrid \M_{u_2})(\X)
=
\frac{1-q^m}{1-q^{n+m}}\,\,
\frac{q^n-t^{\Nmax(u_1)}}{1-t^{\Nmax(u_1)}}\,\,
\M_{u_1}(\X) \M_{u_2}(\X).
\end{equation}
\end{theorem}

\Proof
Thanks to the combinatorial interpretation of $\M_u(\X)$ in terms of signed
words (Proposition~\ref{Mu-col} and Lemma~\ref{lem-1mq}), one only has to
check what happens to the major index of the evaluation of signed words in
the cases of the left, middle, or right tridendriform products.
The analysis is  similar to that done for $\FQSym$ in the previous sections.
% but actually much simpler than this one.
\qed

\begin{example}{\rm
Note that the left tridendriform product does not depend on the actual values
of $w_1$ but only on its length. Indeed, one can check that
\begin{equation}
(\M_{111} \gautrid \M_{2122})(\X)
  = \frac{1-q^3}{1-q^7}\,\, \frac{q^4-t^3}{1-t^3}\,
\M_{111}(\X) \M_{2122}(\X)
\end{equation}
\begin{equation}
(\M_{132} \gautrid \M_{2122})(\X)
  = \frac{1-q^3}{1-q^7}\,\, \frac{q^4-t^3}{1-t^3}\,
\M_{132}(\X) \M_{2122}(\X) 
\end{equation}
But the result depends on the number of maximum of $w_2$:
\begin{equation}
(\M_{121} \gautrid \M_{3122})(\X)
  = \frac{1-q^3}{1-q^7}\,\, \frac{q^4-t}{1-t}\,
\M_{121}(\X) \M_{3122}(\X)
\end{equation}
One can check on these examples the relation
of dendriform trialgebras:
$\M_u\M_v = \M_u\gautrid \M_v + \M_u\miltrid\M_v + \M_u\droittrid\M_v$:
\begin{equation}
(\M_{1212} \gautrid \M_{33231})(\X)
  = \frac{1-q^4}{1-q^9}\,\, \frac{q^5-t^3}{1-t^3}\,
\M_{1212}(\X) \M_{33231}(\X)
\end{equation}
\begin{equation}
(\M_{1212} \miltrid \M_{33231})(\X)
  = \frac{(1-q^4)(1-q^5)}{1-q^9}\,\, \frac{1-t^5}{(1-t^2)(1-t^3)}\,
\M_{1212}(\X) \M_{33231}(\X)
\end{equation}
\begin{equation}
(\M_{1212} \droittrid \M_{33231})(\X)
  = \frac{1-q^5}{1-q^9}\,\, \frac{q^4-t^2}{1-t^2}\,
\M_{1212}(\X) \M_{33231}(\X).
\end{equation}
}
\end{example}

%%%%%%%%%%%%%%%%%%%%%%%%%%%%%%%%%%%%%%%%%%%%%%%%%%%%%%%%%%%%%%%%%%%%%%%%%%%%%%%
%%%%%%%%%%%%%%%%%%%%%%%%%%%%%%%%%%%%%%%%%%%%%%%%%%%%%%%%%%%%%%%%%%%%%%%%%%%%%%%
%%%%%%%%%%%%%%%%%%%%%%%%%%%%%%%%%%%%%%%%%%%%%%%%%%%%%%%%%%%%%%%%%%%%%%%%%%%%%%%
\section{The free dendriform trialgebra}

%%%%%%%%%%%%%%%%%%%%%%%%%%%%%%%%%%%%%%%%%%%%%%%%%%%%%%%%%%%%%%%%%%%%%%%%%%%%%%%
\subsection{A subalgebra of $\WQSym$}

Recall that $\TD$, the Loday-Ronco algebra of plane trees~\cite{LRtri},
is naturally a subalgebra of $\WQSym$~\cite{NT06}, the embedding being
\begin{equation}
\MM_T(A)=\sum_{\TT(u)=T}\M_u(A)\,,
\end{equation}
where $\TT(u)$ is the decreasing plane tree associated with $u$~\cite{NT06}.
Hence, $\MM_T(\X)$ is well-defined.

$\TD$ was originally defined~\cite{LRtri} as the free tridendriform algebra on
one generator as follows: if $T$ is a planar tree and $T_1$, $\dots$, $T_k$
are its subtrees, then
\begin{equation}
\label{defTD}
\MM_T = (\MM_{T_1} \droittrid \MM_1 \gautrid \MM_{T_2}) \miltrid
(\MM_1\gautrid\MM_{T_3}) \miltrid\dots 
(\MM_1\gautrid\MM_{T_k}).
\end{equation}

%%%%%%%%%%%%%%%%%%%%%%%%%%%%%%%%%%%%%%%%%%%%%%%%%%%%%%%%%%%%%%%%%%%%%%%%%%%%%%%
\subsection{$(q,t)$-hooks}

Let $T$ be a plane tree.
Let $\IN(T)$ denote all internal nodes of $T$ except the root.
Let us define a \emph{region} of $T$ as any part of the plane between two
edges coming from the same vertex. The regions are the places where one writes
the values of a packed word when inserting it (see~\cite{NT06}).
For example, with $w=243411$, one gets
\begin{equation}
\vcenter{\xymatrix@C=0.5mm@R=4mm{
*{} & *{} & *{} & *{} & {}\ar@{-}[dlll]\ar@{-}[d]\ar@{-}[drrrr] \\
*{} &  {}\ar@{-}[dl]\ar@{-}[dr] & *{}
& *{4} & {}\ar@{-}[dl]\ar@{-}[dr] & *{}
& *{4} & *{} & {}\ar@{-}[dll]\ar@{-}[d]\ar@{-}[drr] \\ 
 {}
& *{2} &  {} &  {}
& *{3} &  {} &  {}
& *{1} &  {}
& *{1} &  {} \\
}}
\end{equation}

\begin{theorem}
Let $T$ be a plane tree with $n$ regions. Then
\begin{equation}
\MM_T(\X) =
\frac{1-t^{a(r)-1}}{1-q^n}
\prod_{i\in\IN(T)} \frac{q^{r(i)} - t^{a(i)-1} }{1-q^{s(i)-1}},
\end{equation}
where $a(i)$ is the arity of $i$ and $r(i)$ the number of regions of
$T$ below $i$.
\end{theorem}

\Proof
This is obtained by applying the  tridendriform operations  in
$\WQSym$, thanks to~(\ref{defTD}).
\qed

Writing for each node the numerator of its $(q,t)$ contribution, one has, for
example:
\begin{equation}
\vcenter{\xymatrix@C=0.5mm@R=4mm{
*{} & *{} & *{} & *{} & {1-t^2}\ar@{-}[dlll]\ar@{-}[d]\ar@{-}[drrrr] \\
*{} &  {q-t}\ar@{-}[dl]\ar@{-}[dr] & *{}
& *{} & {q-t}\ar@{-}[dl]\ar@{-}[dr] & *{}
& *{} & *{} & {q^2-t^2}\ar@{-}[dll]\ar@{-}[d]\ar@{-}[drr] \\ 
 {}
& *{} &  {} &  {}
& *{} &  {} &  {}
& *{} &  {}
& *{} &  {} \\
}}
\end{equation}
\begin{equation}
\vcenter{\xymatrix@C=0.5mm@R=4mm{
*{} & *{} & *{} & *{} & {1-t}\ar@{-}[dlll]\ar@{-}[drrrr] \\
*{} &  {q-t}\ar@{-}[dl]\ar@{-}[dr] & *{}
& *{} & *{} & *{}
& *{} & *{} & {q^4-t}\ar@{-}[dll]\ar@{-}[drr] \\ 
 {} & *{} &  {} &  {} & *{} &  {} &  {q^3-t^3}\ar@{-}[dll]\ar@{-}[dl]%
\ar@{-}[dr] \ar@{-}[drr]
& *{} &  {} & *{} &  {} \\
 {} & *{} &  {} &  {} & *{} &  {} & {}
& *{} &  {} & *{} &  {} \\
}}
\end{equation}

%%%%%%%%%%%%%%%%%%%%%%%%%%%%%%%%%%%%%%%%%%%%%%%%%%%%%%%%%%%%%%%%%%%%%%%%%%%%%%%
%%%%%%%%%%%%%%%%%%%%%%%%%%%%%%%%%%%%%%%%%%%%%%%%%%%%%%%%%%%%%%%%%%%%%%%%%%%%%%%
%%%%%%%%%%%%%%%%%%%%%%%%%%%%%%%%%%%%%%%%%%%%%%%%%%%%%%%%%%%%%%%%%%%%%%%%%%%%%%%
%%%%%%%%%%%%%%%%%%%%%%%%%%%%%%%%%%%%%%%%%%%%%%%%%%%%%%%%%%%%%%%%%%%%%%%%%%%%%%%

\end{document}